\documentclass[12pt,reqno]{amsart}
\usepackage{amsthm,amsfonts,amssymb,amscd}
\usepackage[all]{xy}

\newtheorem{theorem}[subsection]{Theorem}

\newtheorem{proposition}[subsection]{Proposition}

\newtheorem{lemma}[subsection]{Lemma}
\newtheorem{corollary}[subsection]{Corollary}

\newtheorem{lemma1}[subsubsection]{Lemma}
\newtheorem{prop1}[subsubsection]{Proposition}

\theoremstyle{definition}

\newtheorem{proposition-definition}[subsection]{Proposition-Definition}

\theoremstyle{remark}
\newtheorem{remark}[subsection]{Remark}

\newcommand{\fatdot}{{\scriptscriptstyle \bullet}}
\newcommand{\ccup}{\,{\scriptstyle \cup}\,}
\newcommand{\cccup}{\,{\scriptscriptstyle \cup}\,}
\newcommand{\ccap}{\,{\scriptstyle \cap}\,}
\newcommand{\dual}{{{\scriptscriptstyle \vee}}}

\newcommand{\At}{\operatorname{At}\nolimits}

\newcommand{\codim}{\operatorname{codim}\nolimits}

\newcommand{\EXT}{{\:\mathcal Ext\:}}
\newcommand{\Ext}{\operatorname{Ext}\nolimits}
    \newcommand{\Tor}{\operatorname{Tor}\nolimits}
    \newcommand{\Coh}{\operatorname{Coh}\nolimits}
    \newcommand{\Vect}{\operatorname{Vect}\nolimits}

\newcommand\Hilb{{\operatorname{Hilb}\nolimits}}

\newcommand{\RHOM}{{{\mathcal R}{\mathcal H}om\:}}
\newcommand{\HOM}{{{\mathcal H}om\:}}
\newcommand{\Hom}{\operatorname{Hom}\nolimits}
\newcommand{\RHom}{\operatorname{RHom}\nolimits}
\newcommand{\RGamma}{\operatorname{R\Gamma}\nolimits}
\newcommand{\id}{\operatorname{id}\nolimits}

\newcommand\ob{{\operatorname{ob}\nolimits}}
\newcommand{\Pf}{\operatorname{Pf}\nolimits}

\newcommand{\Pic}{\operatorname{Pic}\nolimits}
\newcommand{\pr}{\operatorname{pr}\nolimits}

\newcommand\sm{{\operatorname{sm}\nolimits}}
\newcommand{\Spec}{\operatorname{Spec}}

\newcommand{\Tr}{\operatorname{Tr}}

\newcommand{\CC}{{\mathbb C}}

\newcommand{\ZZ}{{\mathbb Z}}

\newcommand{\PP}{{\mathbb P}}

    \newcommand{\LL}{{\mathbb L}}

\newcommand{\fc}{{\mathfrak c}}
\newcommand{\FM}{{\mathfrak M}}

\newcommand{\AAA}{{\mathcal A}}
\newcommand{\OOO}{{\mathcal O}}
    \newcommand{\CO}{{\mathcal O}}

\newcommand{\III}{{\mathcal I}}
    \newcommand{\CI}{{\mathcal I}}
\newcommand{\GGG}{{\mathcal G}}
    \newcommand{\CG}{{\mathcal G}}
\newcommand{\EEE}{{\mathcal E}}
    \newcommand{\CE}{{\mathcal E}}
    \newcommand{\D}{{\mathcal D}}

    \newcommand{\CH}{{\mathcal H}}

\newcommand{\FFF}{{\mathcal F}}
\newcommand{\JJJ}{{\mathcal J}}
    \newcommand{\CF}{{\mathcal F}}
\newcommand{\CCC}{{\mathcal C}}
\newcommand{\KKK}{{\mathcal K}}
\newcommand{\NNN}{{\mathcal N}}
    \newcommand{\CN}{{\mathcal N}}

\newcommand{\PPP}{{\mathcal P}}
    \newcommand{\CP}{{\mathcal P}}

\newcommand{\TTT}{{\mathcal T}}
    \newcommand{\CT}{{\mathcal T}}

\newcommand\alp{\alpha}

\newcommand\eps{\epsilon}

\newcommand\si{\sigma}

\newcommand\impl{{\ \Longrightarrow\ }}
\newcommand{\into}{\hookrightarrow}

\newlength{\rrrr}

\newcommand\lra{{\longrightarrow}}

\newcommand\rar{\rightarrow}

\renewcommand{\bar}[1]{\overline{#1}}

    \newcommand{\FMY}{{P(Y)}}

    \newcommand{\lotimes}{\mathbin{\mathop{\otimes}\limits^{\mathrm{L}}}}

\author{A. Kuznetsov}

\address{A. K.: Algebra Section, Steklov Mathematical Institute, 8 Gubkin str., Moscow 119991 Russia}
\email{akuznet@mi.ras.ru}

\thanks{A.K. was partially supported by RFFI grant 05-01-01034, INTAS 05-1000008-8118,
Russian Presidential grant for young scientists No. MK-6122.2006.1,
CRDF Award No. RUM1-2661-MO-05, the Russian Science Support Foundation,
and gratefully acknowledges the support of the Pierre Deligne grant based on his 2004 Balzan prize in mathematics.}

\author{D. Markushevich}

\address{D. M.: Math\'ematiques - b\^{a}t. M2, Universit\'e Lille 1,
F-59655 Villeneuve d'Ascq Cedex, France}
\email{markushe@math.univ-lille1.fr}


\title{Symplectic structures on moduli spaces of sheaves\\
via the Atiyah class}

\begin{document}
\begin{abstract}
It is proven that the composition of the Yoneda coupling
with the semiregularity map is a closed 2-form
on moduli spaces of sheaves. Two examples are given when
this 2-form is symplectic. Both of them are moduli spaces of
torsion sheaves on the cubic 4-fold $Y$. The first example
is the Fano scheme of lines in $Y$. Beauville and Donagi
showed that it is symplectic but did not construct
an explicit symplectic form on it. We prove that our construction
provides a symplectic form. The other example is
the moduli space of torsion sheaves which
are supported on the hyperplane sections $H\cap Y$ of $Y$ and
are cokernels of the Pfaffian
representations of $H\cap Y$.
\end{abstract}
\maketitle

\section*{Introduction}

The existence of
symplectic or Poisson structures on various
moduli spaces in differential and algebraic geometry
is rather a general phenomenon. Starting from the seminal
work of Atiyah--Bott \cite{AB}, such structures appeared
mostly via the method of symplectic reduction, which was
further developed in order to produce K\"ahler
and hyper-K\"ahler structures \cite{KN}. Nowadays, many
results on holomorphic Poisson
structures have been obtained by techniques of algebraic geometry
for moduli spaces of sheaves \cite{Muk-1}, \cite{Bot-1}, \cite{O'G},
\cite{K}
or related objects, including Higgs pairs \cite{Hi}, \cite{Bot-2},
\cite{BR}, parabolic bundles \cite{Bot-3}, regular or meromorphic connections \cite{IIS}.

Mukai \cite{Muk-1} proved that any moduli space of simple sheaves on
a K3 or abelian surface has a nondegenerate holomorphic 2-form.
Its closedness was proved later in \cite{Muk-2},
\cite{O'G},  \cite{Ran} for vector bundles, in \cite{Bot-1} for torsion-free sheaves,
and in \cite{HL} for arbitrary sheaves. Mukai's result
was extended in \cite{Tyu} to moduli spaces of vector bundles over surfaces of general type
and over Poisson surfaces; a more thorough study of the Poisson case
was accomplished in \cite{Bot-1}.
Higher-dimensional generalizations were obtained in
\cite{K} and \cite{Ran}. Kobayashi proved that moduli spaces of simple vector bundles on a
hyper-K\"ahler manifold are holomorphically symplectic. Ran extended the closedness
result for the 2-form of Mukai's type to simple vector bundles on compact
complex manifolds.
In all these situations, the 2-form or a Poisson bivector on the moduli space of sheaves
is induced by that on the base space of the sheaves.

Beauville--Donagi \cite{BD} discovered that the variety $F(Y)$
of lines in a smooth cubic 4-fold $Y\subset\PP^5$ is holomorphically
symplectic. Their proof is indirect: they identified $F(Y)$
for a special $Y$ with the length-2 punctual Hilbert scheme
of a K3 surface, which is known to be
an irreducible symplectic manifold. This means that
it is compact, simply connected, hyper-K\"ahler and has a
unique holomorphic symplectic structure. Then the assertion for
any smooth 4-dimensional cubic follows by a deformation argument:
any K\"ahler deformation of an irreducible symplectic
manifold is also irreducible symplectic. This approach does
not provide a recipe for {\em constructing} a symplectic form on
$F(Y)$.

One can interprete $F(Y)$ as the moduli space
parameterizing the structure sheaves of all the lines on $Y$.
The authors of \cite{MT2} found another holomorphically symplectic
moduli space of sheaves $\FMY$ on the cubic 4-fold $Y$.
It parametrizes
the torsion sheaves which are rank-2 vector bundles on the
hyperplane sections of $Y$ with Chern numbers $c_1=0$, $c_2=2$
(see Theorem~\ref{bundles-e} for other equivalent characterizations
of these sheaves).
These examples both differ from the previous ones by the fact
that the symplectic structures on them are nomore
induced by a symplectic structure on $Y$ itself: $Y$ has no
holomorphic forms at all.

One of the objectives of this paper is to find a general construction
of closed $p$-forms on moduli spaces which yields the symplectic structures in the
above two examples. We produce such a construction, involoving the
Yoneda product with the Atiyah class
of the sheaves and valid over the smooth locus of all the moduli spaces of
sheaves on an arbitrary smooth complex projective variety.

The other objective of this paper is to gather some general techniques for working with the Atiyah class of non-locally-free sheaves that might be useful in the study of our $p$-forms
in other examples. In particular, we relate the Atiyah class of a torsion sheaf
to its linkage class (see Sect. 3). The linkage class is crucial in our proof
of the nondegeneracy of the 2-form on $P(Y)$. This tool seems to be
not widely known; it represents a certain novelty. As applications
of the general techniques, we also give shortcut
formulas for 2-forms on the Hilbert scheme of l.~c.~i. subschemes of
a given smooth projective variety (Sect. 6) and an explicit computation
of the Beauville--Donagi symplectic form on $F(Y)$ (see (\ref{2-form-in-coo}),
(\ref{2-form-in-coo2})).

Two other
approaches to an explicit formula for the Beauville--Donagi
form are described in \cite{dJS}, \cite{IMan}. The authors of \cite{dJS} construct
2-forms on the moduli spaces of degree-$d$ rational curves in the cubic 4-fold.
Their construction is specific for a cubic 4-fold and does not involve main technical
tools of our work, neither the Atiyah class, nor the linkage class. Their moduli
spaces contain open sets that can be identified as moduli of sheaves, parametrizing
the structure sheaves of the rational curves. The shortcut formula for Hilbert
schemes from Sect. 6 implies that our 2-form is proportional to that of
de Jong--Starr over these open sets. The approach of \cite{IMan} is completely different
and uses the embedding $F(Y)\subset G(2,6)$.

Now we will briefly describe our construction of $p$-forms for the case $p=2$.
It involves the following steps.
The tangent space to the moduli space at a point $[\FFF]$
representing a stable (or just simple) sheaf $\FFF$ is canonically
isomorphic to $\Ext^1(\FFF,\FFF)$, so we have to associate a complex number
to two elements of $\Ext^1(\FFF,\FFF)$. The first step is the Yoneda
coupling
$$
\Ext^1(\FFF,\FFF)\times\Ext^1(\FFF,\FFF)\lra\Ext^2(\FFF,\FFF).
$$
This bilinear map is skew-symmetric whenever $[\FFF]$ is a nonsingular
point of the moduli space. Indeed, according to \cite{Ar} and \cite{Muk-1}, the quadratic term of the obstruction map $\Ext^1(\FFF,\FFF)\rar\Ext^2(\FFF,\FFF)$ is the Yoneda square,
so, if $[\FFF]$ is smooth, then the Yoneda square vanishes and the Yoneda
coupling is skew-symmetric.

When the base space of the sheaves is a symplectic surface $S$
with a symplectic form $\omega^{2,0}\in H^0(S,\Omega^2_S)$,
Mukai composes the Yoneda coupling with the map
$$
\Ext^2(\FFF,\FFF)\xrightarrow{\Tr}H^2(S,\OOO_S)
\xrightarrow{\ccup \omega^{2,0}}H^2(S,\Omega^2_S)=\CC,
$$
and this ends the construction in the surface case. Over a $n$-dimensional
base $Y$, we insert an intermediate step: compose the Yoneda
coupling with an exterior
power of the Atiyah class $\At (\FFF)\in \Ext^1(\FFF,\FFF\otimes\Omega^1_Y)$:
\begin{equation}\label{step2}
\Ext^2(\FFF,\FFF)\ni\xi\mapsto {\At(\FFF)^{\wedge q}}\circ\xi,\ \
\Ext^q(\FFF,\FFF\otimes\Omega^q_Y)\times\Ext^2(\FFF,\FFF)
\xrightarrow{\ \circ\ } \Ext^{q+2}(\FFF,\FFF\otimes\Omega^q_Y).
\end{equation}
The exponent $q$ should be chosen in such a way that $h^{q,q+2}(Y)\neq 0$.
Then we pick up an element $\omega=\omega^{n-q,n-q-2}\in H^{n-q-2}(\Omega_Y^{n-q})$, and, to
end up in $\CC$, compose with the map
\begin{equation}\label{step3}
\Ext^{q+2}(\FFF,\FFF\otimes\Omega^q_Y)\xrightarrow{\Tr}
H^{q+2}(\Omega_Y^q)\xrightarrow{\ccup \omega^{n-q,n-q-2}}
H^{n}(\Omega_Y^{n})=\CC,
\end{equation}
which provides the 2-form $\alp_\omega$ on the moduli space.
The idea to couple $\Ext^2(\FFF,\FFF)$ with powers of $\At(\FFF)$
was used by Buchweitz--Flenner \cite{BuF1}, \cite{BuF2} to define
an analog of the Bloch semiregularity map for deformations of sheaves.
Thus our approach combines the ideas of Mukai and Buchweitz--Flenner.

Unlike Mukai's case, in which the nondegeneracy of the 2-form immediately follows
from the Serre duality, our forms may be degenerate, and
it is not so easy to prove that they are nondegenerate or even nonzero
in particular examples.
For a cubic 4-fold, $h^{1,3}=1$, so that our construction provides
a unique, up to a constant factor, 2-form $\alp$ on every moduli space
of sheaves on $Y$. We prove the following general sufficient
condition: $\alp$ is nondegenerate at $[\FFF]$ if
\begin{equation}\label{criterion}
H^i(\FFF)=H^i(\FFF(-1))=H^i(\FFF(-2))=0 \ \ \mbox{for all}\ \  i\in\ZZ.
\end{equation}
This condition is verified for all sheaves in $\FMY$, but not
for the structure sheaves $\OOO_\ell$ of lines $\ell\subset Y$.
However, we manage to apply this criterion to $F(Y)$ upon replacing
$\OOO_\ell$ by the second syzygy sheaf of $\OOO_\ell(1)$,
see the last part of Section~\ref{sc4} and Section~\ref{sdim4} for more details.

The adequate techniques for the proof of the nondegeneracy criterion
are those of derived category. The real reason of its validity
is that the full triangulated subcategory $\CCC_Y\subset\D^b(Y)$
of those $\FFF$
satisfying (\ref{criterion}) is a kind of deformation of the
derived category of a K3 surface (see \cite{Ku2}), and the moduli of sheaves
in it behave like moduli of sheaves on a K3 surface. The K3-type Serre
duality $\Ext^i(\FFF,\GGG)\cong \Ext^{2-i}(\GGG,\FFF)^\dual$ on $\CCC_Y$ is defined
via the composition with the linkage class
$\eps_\GGG\in\Ext^2(\GGG,\GGG\otimes\NNN^\dual_{Y/\PP^5})$:
$$
\Ext^i(\FFF,\GGG)\times \Ext^{2-i}(\GGG,\FFF)\xrightarrow{\ \circ\ }
\Ext^{2}(\GGG,\GGG)\xrightarrow{\eps_\GGG\circ\cdot\ }\Ext^{4}(\GGG,\GGG\otimes\Omega_Y^4)\xrightarrow{\Tr}H^4(\Omega_Y^4)=\CC,
$$
where we use an isomorphism $\NNN^\dual_{Y/\PP^5}\cong \Omega_Y^4$.
This implies the nondegeneracy of the forms $\alp$ by virtue of
Theorem \ref{epsat}, which affirms that $\eps_\GGG$ factors through
$\At(\GGG)$.

Remark that both $\CCC_Y$ and $P(Y)$ are quite intriguing objects
associated to any smooth cubic 4-fold $Y$. For special $Y$'s,
they are related to certain K3 surfaces.
We state two conjectures: 1)~$Y$~is rational
if and only if $\CCC_Y$ is equivalent to the derived category of
a K3 surface, and 2)~$P(Y)$~has a compactification, deformation
equivalent to the O'Grady irreducible symplectic 10-fold \cite{OG}.
See \cite{Ku2} for some evidence confirming the first conjecture.
The second conjecture would follow from proving that when
$Y$ is a Pfaffian cubic, then the integral Fourier--Mukai
functor with kernel $\EEE$ described in loc. cit., Sect. 8, maps
$P(Y)$ to the O'Grady 10-fold associated to the K3 surface $S$,
projectively dual to $Y$.

The construction of  $p$-forms for any $p$ is obtained
by taking in its first step the $p$-linear Yoneda
product on $\Ext^1(\FFF,\FFF)$ in place of the bilinear one.
In the sequel, we limit ourselves to $p=2$ for notational convenience,
and also in view of the particular importance of this case,
including symplectic structures. We will only briefly comment on
what is known for other $p$'s. If $p=1$,
our map $\omega\mapsto\alp_\omega$ is the
adjoint of the infinitesimal Abel-Jacobi map of Buchweitz--Flenner, a generalization
of Griffiths' map for a Hilbert scheme as defined in \cite{Gri}, Theorem 2.5.
The formula of Griffiths
is the exact $(p=1)$-analog of our formula (\ref{beta}).
It is also worth noting that Welters in \cite{W}, 2.8, in order
to calculate the infinitesimal Abel-Jacobi map for curves
on a 3-fold $X$, embeds $X$ into a 4-fold $W$ and obtains
a description, equivalent to our recipe of taking product with the linkage
class of the embedding $X\subset W$.

For $p\geq 3$, we refer to \cite{T}, where a cubic Yoneda coupling is used
to produce a 3-form on a particular moduli space of sheaves on a Calabi--Yau 3-fold $X$
fibered in $K3$ surfaces $S_t$,
which makes this moduli space into another Calabi-Yau 3-fold $\hat X$, a kind of mirror
of $X$. It is interesting to
note that the proof of the nondegeneracy of this 3-form in \cite{T} implicitly
uses the linkage class of the embedding $S_t\subset X$, appearing as the
connecting homomorphism $\delta$ in Lemma 3.42.

Now we will describe the contents of the paper by sections.
In Section 1, we collect reminders on the tools
needed in the sequel: trace map, functors $Li^*$, $i^!$
and duality for a closed embedding $i:Z\into Y$,
evaluation (or integral transform) and the
Atiyah class. In Section 2, we provide the construction of the
2-form $\alp$ and prove that it is closed in adapting to our case
the proof of the closedness of Mukai's form in \cite{HL}.
In Section 3, we define the linkage class $\eps_\FFF$ of a sheaf
$\FFF$ supported on a locally complete intersection subscheme
$Y$ in a variety $M$ and
show that $\eps_\FFF$ factors through
$\At(\FFF)$. In Section 4, we show that on a cubic 4-fold,
the product with
$\eps_\GGG$ induces an isomorphism from $\Ext^{\fatdot}(\FFF,\GGG)$
onto $\Ext^{\fatdot+2}(\FFF,\GGG(-3))$ whenever $\FFF,\GGG\in\CCC_Y$,
which implies the nondegeneracy of the 2-form $\alp$ on any moduli
space parametrizing sheaves from $\CCC_Y$. Section 5 represents
the family $F(Y)$ of lines on a cubic 3-fold as a connected
component of a moduli space of sheaves from $\CCC_Y$, hence the nondegeneracy
of the 2-form $\alp$ on it is a consequence of the results
of the previous section. Section 6 provides a simpler shortcut formula for calculation
of the 2-form in the case when the moduli space under consideration
is (the partial compactification of) the relative Picard
of some Hilbert scheme of equidimensional l.~c.~i. subschemes
of $Y$. This formula implies that the 2-form lifts from the
Hilbert scheme of $Y$. As an application, we deduce explicit
formulas in coordinates for the case of lines in a cubic 4-fold.
The concluding Section 7 describes the 10-dimensional moduli space $\FMY$.

\bigskip

{\sc Acknowledgements.} The second author thanks Hubert Flenner for
discussions. He also acknowledges the hospitality of
the Max-Planck-Institut f\"ur Mathematik in Bonn and the Mittag-Leffler
Institute in Stockholm, where was done the work on the present paper.

\section{Preliminaries}

\subsection{Notation and conventions}

Throughout the paper we use the field of complex numbers $\CC$ as the base field.
Certainly, all results remain true for any algebraically closed field of zero characteristic.
On the other hand, the Hodge theory is used in the proof of closedness of the forms,
so this probably fails in positive characteristic. By an algebraic variety we mean
an integral separated scheme of finite type over the base field.

Given an algebraic variety $Y$ we denote by $\Coh(Y)$ the abelian category of coherent sheaves on $Y$,
and by $\D^b(\Coh(Y))$ its bounded derived category. It is defined (see \cite{Ve}) as the localization
of the homotopy category of bounded complexes of coherent sheaves with respect to the class of quasiisomorphisms
of complexes. There are also some unbounded versions of the derived category:
the bounded above derived category $\D^-(\Coh(Y))$,
the bounded below derived category $\D^+(\Coh(Y))$,
and the unbounded derived category $\D  (\Coh(Y))$.
The derived category is triangulated, it comes equipped with the {\em shift functor}\/ $F \mapsto F[1]$
(induced by the shift of grading on complexes) and with a class of {\em distinguished}\/ (or {\em exact})
{\em triangles}, sequences of the form $\xymatrix@1{F_1 \ar[r] & F_2 \ar[r] & F_3 \ar[r] & F_1[1]}$
(generalizing short exact sequences of complexes), satisfying a number of axioms.

All the standard functors on categories of coherent sheaves give rise to their derived functors
between the derived categories (see, e.g.~\cite{Ha-1}). Derived functors are compatible with triangulated structures,
i.e. commute with the shift functor and take exact triangles to exact triangles. In particular,
the tensor product $\otimes:\Coh(Y) \times\Coh(Y) \to \Coh(Y)$ gives rise
to the derived tensor product $\lotimes:\D^-(\Coh(Y))\times\D^-(\Coh(Y)) \to \D^-(\Coh(Y))$,
the functor of local homomorphisms $\HOM:\Coh(Y)^\circ \times\Coh(Y) \to \Coh(Y)$ gives rise
to the derived local-$\Hom$ functor $\RHOM:\D^-(\Coh(Y))^\circ\times\D^+(\Coh(Y)) \to \D^+(\Coh(Y))$,
and the functor of global homomorphisms $\Hom:\Coh(Y)^\circ \times\Coh(Y) \to \Vect$ gives rise
to the derived global $\Hom$-functor $\RHom:\D^-(\Coh(Y))^\circ\times\D^+(\Coh(Y)) \to \D^+(\Vect)$,
where $\Vect$ stands for the category of vector spaces. Similarly, given a proper map $f:X \to Y$
of algebraic varieties,
the pushforward functor $f_*:\Coh(X) \to \Coh(Y)$ gives rise to
the derived pushforward $Rf_*:\D^b(\Coh(X)) \to \D^b(\Coh(Y))$ and
the pullback functor $f^*:\Coh(Y) \to \Coh(X)$ gives rise to
the derived pullback $Lf^*:\D^-(\Coh(Y)) \to \D^-(\Coh(X))$.
In particular, when $f:X \to \Spec\CC$ is the projection to the point,
the pushforward functor $f_*$ is nothing but the functor of global sections $\Gamma(X,-)$
and its derived functor is denoted by $\RGamma(X,-)$.

Given an object $F \in \D(\Coh(Y))$ we can consider its $k$-th cohomology sheaf $\CH^k(F)$.
The cohomology sheaves of the derived functors applied to sheaves are the classical derived functors,
e.g. $\CH^{-k}(F\lotimes G) = \Tor_k(F,G)$, $\CH^k(\RHOM(F,G)) = \EXT^k(F,G)$, $\CH^k(\RHom(F,G)) = \Ext^k(F,G)$,
$\CH^k(Rf_*(F)) = R^kf_*(F)$, $\CH^{-k}(Lf^*(G)) = L_kf^*(G)$, and $\CH^k(\RGamma(X,F)) = H^k(X,F)$.
This gives a useful interpretation of $\Ext$'s.
We have
$$\Ext^k(F,G) = \CH^k(\RHom(F,G)) = \CH^0(\RHom(F,G)[k]) = \CH^0(\RHom(F,G[k])) = \Hom(F,G[k]),$$
so we can consider an element of the space $\Ext^k(F,G)$ as a morphism $F \to G[k]$.
In this interpretation the Yoneda multiplication of $\Ext$'s corresponds to the composition of $\Hom$'s.

The standard isomorphisms between functors give rise to isomorphisms between derived functors.
E.g., $\Hom(F,G) \cong \Gamma(X,\HOM(F,G))$ gives $\RHom(F,G) \cong \RGamma(X,\RHOM(F,G))$.
Considering the cohomology sheaves, these isomorphisms give rise to spectral sequences,
such as the local-to-global spectral sequence
\begin{equation}\label{ltg}
E_2^{p,q}=H^p(X,\EXT^q(\FFF,\GGG))\impl \Ext^{p+q}(\FFF,\GGG).
\end{equation}

Since in general the derived functors may take an object of the bounded derived category
to an unbounded object, the conditions ensuring boundedness are very useful.
We will mention two of them. An object $F \in \D^b(\Coh(Y))$ is called a {\em perfect complex}\/
if it is locally quasiisomorphic to a finite complex of locally free sheaves of finite rank.
If $F$ is a perfect complex then the functors $\RHOM(F,-)$ and $\RHom(F,-)$ preserve boundedness.
Moreover, the derived pullback of a perfect complex is a perfect complex.
Perfect complexes form a triangulated subcategory of $\D^b(\Coh(Y))$ called the category
of perfect complexes.

A map $f:X \to Y$ is said to have finite $\Tor$-dimension if the structure sheaf $\CO_X$
has finite $\Tor$-dimension over $\CO_Y$, i.e. for any point $x\in X$ the local ring
$\CO_{X,x}$ admits a finite flat resolution over $\CO_{Y,f(x)}$.
If $f:X \to Y$ has finite $\Tor$-dimension then the derived pullback functor preserves boundedness
and the derived pushforward functor takes perfect complexes to perfect complexes.

\subsection{Traces}

Let $Y$ be an algebraic variety and $E$ a vector bundle on $Y$.
For every coherent sheaf $\CF$ on $Y$ which is a perfect complex
consider the composition
$$
\RHOM(\CF,\CF\otimes E) \cong \CF^\dual\lotimes\CF\otimes E \to E,
$$
where $\CF^\dual = \RHOM(\CF,\CO_Y)$ is the derived dual of $\CF$.
The first map above is the canonical isomorphism (it uses perfectness of $\CF$)
and the second is the ``contraction'' map.
Taking the $k$-th cohomology, we obtain a natural map
$$
\Tr:\Ext^k(\CF,\CF\otimes E) \to H^k(Y,E),
$$
the {\em trace map} (see~\cite{Ill}).

The most important property of the trace is {\em additivity}:
if $\xymatrix@1@C=15pt{\CF_1 \ar[r]^{\phi_1} & \CF_2 \ar[r]^{\phi_1} & \CF_3 \ar[r]^{\phi_3} & \CF_1[1]}$
is a distinguished triangle and a collection $\mu_i \in \Ext^k(\CF_i,\CF_i\otimes E)$
is compatible with the triangle (i.e. the diagram
$$
\xymatrix@C=35pt{
\CF_1 \ar[r]^{\phi_1} \ar[d]_{\mu_1} & \CF_2 \ar[r]^{\phi_2} \ar[d]_{\mu_2} & \CF_3 \ar[r]^{\phi_3} \ar[d]_{\mu_3} & \CF_1[1] \ar[d]_{-\mu_1} \\
\CF_1\otimes E[k] \ar[r]^{\phi_1\otimes 1_E} & \CF_2\otimes E[k] \ar[r]^{\phi_2\otimes 1_E} & \CF_3\otimes E[k] \ar[r]^{\phi_3\otimes 1_E} & \CF_1\otimes E[k+1]
}
$$
is commutative), then
$$
\Tr(\mu_1) - \Tr(\mu_2) + \Tr(\mu_3) = 0
$$
in $H^k(Y,E)$.

Another important property is {\em multiplicativity}:
if $\mu \in \Ext^k(\CF,\CF\otimes E)$ and $\varphi \in \Ext^l(E,E')$ then
$$
\varphi\circ\Tr(\mu) = \Tr((\id_\CF\otimes\varphi)\circ\mu)
$$
in $H^{k+l}(Y,E')$.

\subsection{Sheaves on a subvariety}

Let $i:Z \into Y$ be a closed embedding.
If $Z \subset Y$ is a locally complete intersection,
we denote by $\CN_{Z/Y}$ the normal bundle of $Z$ in $Y$.
Now let us compute the cohomology sheaves $L_ki^*\CF$ of the derived pullback functor
for $\CF = i_*F$, where $F$ is a coherent sheaf on $Z$.

\begin{lemma1}\label{pbpf}
If $Z \subset Y$ is a locally complete intersection of codimension $m$ then
we have $L_ki^*i_*F \cong F\otimes\wedge^k\CN^\dual_{Z/Y}$ for $0 \le k\le m$,
and $L_{>m}i^*i_*F = 0$.
\end{lemma1}
\begin{proof}
Since $i$ is a closed embedding, it suffices to check that
$$
i_*L_ki^*i_*F \cong i_*(F\otimes\wedge^k\CN^\dual_{Z/Y}),
\quad\text{for $0\le k\le m$},
\qquad\text{and}\qquad
i_*L_{>m}i^*i_*F = 0.
$$
By the projection formula we have
$i_*L_ki^*i_*F \cong \Tor_k(i_*F,i_*\CO_Z)$.
Since $Z$ is a locally complete intersection,
it can be represented locally as the zero locus of
a regular section of a rank $m$ vector bundle $\CE$ on $Y$.
Therefore, locally we have the Koszul resolution
$$
0 \to \wedge^m\CE^\dual \to \wedge^{m-1}\CE^\dual \to \dots
\to \wedge^2\CE^\dual \to \CE^\dual \to \CO_Y \to i_*\CO_Z \to 0.
$$
Using it to compute $\Tor$-s we see that
$\Tor_k(i_*\CF,i_*\CO_Y) \cong
i_*\CF\otimes\wedge^k\CE^\dual \cong
i_*(\CF\otimes\wedge^k\CE^\dual_{|Z})$
and it remains to note that $\CN^\dual_{Z/Y} \cong \CE^\dual_{|Z}$.
\end{proof}

Now let us compute $\EXT^k(\CF,\CG)$ for $\CF = i_*F$, $\CG = i_*G$,
where $F$ and $G$ are coherent sheaves on $Z$.

\begin{lemma1}\label{extfg}
Assume that $Z \subset Y$ is a locally complete intersection of codimension $m$.
Let $F$ and $G$ be coherent sheaves on $Z$ and assume that $F$ is locally free.
Then we have

$(i)$ $\EXT^k(i_*F,i_*G) \cong
i_*(\wedge^k\CN_{Z/Y}\otimes F^\dual\otimes G)$ if $0\leq k\leq m$, and
$\EXT^k(i_*F,i_*G) = 0$ otherwise.

$(ii)$ Multiplication
$\EXT^l(i_*G,i_*H)\otimes\EXT^k(i_*F,i_*G) \to \EXT^{k+l}(i_*F,i_*H)$
corresponds under isomorphisms $(i)$ to the map
$\wedge^l\CN_{Z/Y}\otimes G^\dual\otimes H\otimes\wedge^k\CN_{Z/Y}\otimes F^\dual\otimes G \to \wedge^{k+l}\CN_{Z/Y}\otimes F^\dual\otimes H$
given by the wedge product $\wedge^l\CN_{Z/Y}\otimes\wedge^k\CN_{Z/Y} \to \wedge^{k+l}\CN_{Z/Y}$ and the contraction
$G^\dual\otimes G \to \CO_Z$.
\end{lemma1}
\begin{proof}
For $(i)$ we note that $R\HOM(i_*F,i_*G) \cong i_*R\HOM(Li^*i_*F,G)$ by the standard adjunction
between the pushforward and the pullback (see~\cite{Ha-1}). On the other hand, by Lemma~\ref{pbpf}
the complex $Li^*i_*F$ has cohomology sheaves
$F\otimes\wedge^k\CN^\dual_{Z/Y}$, which are locally free.
Therefore, $R\HOM(Li^*i_*F,G)$ has cohomology
$(F\otimes\wedge^k\CN^\dual_{Z/Y})^\dual\otimes G \cong \wedge^k\CN_{Z/Y}\otimes F^\dual\otimes G$.

Assertion $(ii)$ is local, so we may assume that $Z$ is the zero locus of
a regular section of a rank $m$ vector bundle $\CE$ on $Y$, and that $F$, $G$ and $H$
are the restrictions from $Y$ to $Z$ of sheaves $\CF$, $\CG$ and $\CH$.
Then the tensor product of $\CF$ and of the Koszul resolution of $i_*\CO_Z$
is a resolution of $i_*F \cong i_*i^*\CF \cong \CF\otimes i_*\CO_Z$,
similarly for $i_*G$, and we use these resolutions to compute $R\HOM$-s.
It is clear that the multiplication $R\HOM(i_*F,i_*G)\otimes R\HOM(i_*G,i_*H) \to R\HOM(i_*F,i_*H)$
is induced by the wedge product $\wedge^k\CE\otimes\wedge^l\CE \to \wedge^{k+l}\CE$
and by the contraction $\CG^\dual\otimes\CG \to \CO_Y$. Restricting
to the cohomology we deduce the claim.
\end{proof}

The local-to-global spectral sequence~\eqref{ltg}
allows us to compute $\Ext^k(\CF,\CG)$ for the sheaves $\CF = i_*F$, $\CG = i_*G$.
For the computation of the Yoneda multiplication on $\Ext$-s the following
lemma is very useful.

\begin{lemma1}\label{ltgm}
The maps $$H^{p_1}(V,\EXT^{q_1}(\CG,\CH)) \otimes H^{p_2}(V,\EXT^{q_2}(\CF,\CG)) \to H^{p_1+p_2}(V,\EXT^{q_1+q_2}(\CF,\CH))$$
induced by the composition on local $\EXT\!$'s and by the cup-product on the cohomology commute with
the differentials of the spectral sequence and
differ from the maps induced by the Yoneda multiplication
$\Ext^{k_1}(\CG,\CH) \otimes \Ext^{k_2}(\CF,\CG) \to \Ext^{k_1+k_2}(\CF,\CH)$
by the sign $(-1)^{p_1q_2}$.
\end{lemma1}
\begin{proof}
Follows from the fact that the isomorphism of functors $\RGamma \circ \RHOM \cong \RHom$
is compatible with multiplication upon an appropriate change of signs
in the double complex on the l.~h.~s.
\end{proof}

\subsection{Traces and duality}\label{traces-and-duality}

Let $f:Z \to Y$ be a projective morphism of finite $\Tor$-dimension.
It is well known that the (derived) pushforward functor
$Rf_*:\D^b(\Coh(Z)) \to \D^b(\Coh(Y))$ has a right adjoint functor,
called the {\em twisted pullback functor}\/ $f^!:\D^b(\Coh(Y)) \to \D^b\Coh(Z))$
(see~\cite{Ha-1}). The twisted pullback functor differs from the usual (derived)
pullback functor by the relative canonical class (and a shift), explicitly
$$
f^!\CF \cong Lf^*\CF \otimes \omega_{Z/Y}[\dim Z - \dim Y].
$$
Since $f^!$ is right adjoint to $Rf_*$ we have a canonical {\em duality} isomorphism
\begin{equation}\label{duality}
\Hom(\CG,f^!\CF) \cong \Hom(Rf_*\CG,\CF)
\end{equation}
for all $\CF\in\D^b(\Coh(Y))$, $\CG \in \D^b(\Coh(Z))$.
In particular, taking $\CG = f^!\CF$ we denote by $\mathsf{can}$
the map $Rf_*f^!\CF \to \CF$ corresponding to the identity map $f^!\CF \to f^!\CF$
under the duality isomorphism.

\begin{lemma1}\label{tp-factors}
Take any morphism $\varphi: \CG \to f^!\CF$. Then the morphism $Rf_*\CG \to \CF$
associated to $\varphi$ by the duality isomorphism coincides with the composition
\begin{equation}\label{factorization}
\xymatrix@1{Rf_*\CG \ar[r]^-{Rf_*\varphi} & Rf_*f^!\CF \ar[r]^-{\mathsf{can}} & \CF}.
\end{equation}
In particular, every map $Rf_*\CG \to \CF$ has a factorization~\eqref{factorization}
for some $\varphi:\CG \to f^!\CF$.
\end{lemma1}
\begin{proof}
This is a standard consequence of adjointness.
\end{proof}

We will use the twisted pullback functor in case when $f$ is a closed embedding
of a locally complete intersection subscheme. So, assume that $Z \subset Y$
is a locally complete intersection of codimension $m$, and $i:Z \to Y$
is the embedding. Then $\omega_{Z/Y} = \wedge^m\CN_{Z/Y}$ and we have
\begin{equation}\label{tpi}
i^!\CF \cong Li^*\CF\otimes\wedge^m\CN_{Z/Y}[-m].
\end{equation}
In particular, if $\CF$ is a vector bundle on $Y$, then $i^!\CF$
is a vector bundle on $Z$ shifted by $-m$.

\begin{prop1}\label{trace_factors}
Let $Y$ be an algebraic variety, $Z \subset Y$ a locally complete intersection subscheme,
$i:Z \to Y$ the embedding, $E$ a vector bundle on $Y$,
and $F$ a coherent sheaf on $Z$ which is a perfect complex.
Then the trace map $\Tr:\Ext^k(i_*F,i_*F \otimes E) \to H^k(Y,E)$ factors through
$H^k(Y,i_*i^!E)\xrightarrow{\mathsf{can}}H^k(Y,E)$.
\end{prop1}

\begin{proof}
By definition, the trace map is induced by the contraction map
$(i_*F)^\dual\lotimes i_*F\otimes E \to E$ and by the projection formula we have
$(i_*F)^\dual\lotimes i_*F\otimes E \cong i_*(Li^*((i_*F)^\dual\otimes E)\lotimes F)$,
hence the contraction map factors through $i_*i^! E \xrightarrow{\mathsf{can}} E$
by Lemma~\ref{tp-factors}. By functoriality of the cohomology,
the trace map factors as well.
\end{proof}

In the case when $k = m$ is the codimension of $Z$ in $Y$, and $F$ is locally free,
the factorization of the trace map can be described rather explicitly.
In this case we can consider the following composition of maps
\begin{multline}\label{e-f}
\Ext^m(i_*F,i_*F\otimes E) \to
H^0(Y,\EXT^m(i_*F,i_*F\otimes E)) \cong \\ \cong
H^0(Z,F^\dual\otimes F\otimes\wedge^m\CN_{Z/Y}\otimes E_{|Z}) \to
H^0(Z,\wedge^m\CN_{Z/Y}\otimes E_{|Z}) \cong \\ \cong
H^m(Z,i^!E) \cong
H^m(Y,i_*i^!E) \xrightarrow{\mathsf{can}}
H^m(Y,E).
\end{multline}
The first map here is the canonical projection,
the second is the isomorphism of Lemma~\ref{extfg}~$(i)$,
the third is the trace map on $Z$,
the fourth is the isomorphism~\eqref{tpi},
the fifth is evident, and the last one is the canonical map.

\begin{prop1}\label{exp-fact}
Let $Y$ be an algebraic variety, $Z \subset Y$ a locally complete intersection subscheme
of codimension $m$, $i:Z \to Y$ the embedding, $E$ a vector bundle on $Y$,
and $F$ a vector bundle on $Z$. Then the composition of the maps in~\eqref{e-f}
coincides with the trace map
$\Tr:\Ext^m(i_*F,i_*F \otimes E) \to H^m(Y,E)$.
\end{prop1}
\begin{proof}
Indeed, by Lemma~\ref{extfg}~$(i)$ the complex
$(i_*F)^\dual\lotimes i_*F\otimes E$ has nontrivial cohomology only in degrees from $0$ to $m$,
while $i_*i^!E$ is concentrated in degree $m$ by~\eqref{tpi}. Therefore, the contraction map
$(i_*F)^\dual\lotimes i_*F\otimes E \to i_*i^!E$ factors through the $m$-th cohomology sheaf
of $(i_*F)^\dual\lotimes i_*F\otimes E$, i.e. through $\EXT^m(i_*F,i_*F\otimes E)$.
The remaining part is evident.
\end{proof}

Finally, we will need the following description of the map
$H^m(Y,i_*i^!E) \xrightarrow{\mathsf{can}} H^m(Y,E)$.

\begin{lemma1}\label{can}
Let $Y$ and $Z$ be smooth varieties, $Z \subset Y$, $m = \codim_Y Z$.
The map $H^m(Y,i_*i^!E) \xrightarrow{\mathsf{can}} H^m(Y,E)$ is dual to the map
\begin{multline*}
H^m(Y,E)^\dual \cong
H^{n-m}(Y,E^\dual\otimes\omega_Y) \to
H^{n-m}(Z,E_{|Z}^\dual\otimes\omega_{Y|Z}) \cong \\ \cong
H^0(Z,E_{|Z}\otimes\omega_{Y|Z}^{-1}\otimes\omega_Z)^\dual \cong
H^0(Z,E_{|Z}\otimes\wedge^m\CN_{Z/Y})^\dual \cong
H^m(Z,i^!E)^\dual \cong
H^m(Y,i_*i^!E)^\dual
\end{multline*}
where the first and the third isomorphisms are given by the Serre duality
on $Y$ and $Z$ respectively, the second map is the restriction from $Y$ to $Z$,
the fourth isomorphism is the adjunction formula for $\omega_Z$,
the fifth is~\eqref{tpi}, and the last one is obvious.
\end{lemma1}
\begin{proof}
The standard relation between the right adjoint and the left adjoint functors.
When the Serre duality takes place, the right adjoint functor (and its canonical map)
is obtained from the left adjoint functor by conjugation with the Serre functors.
\end{proof}

\subsection{Evaluation}

Let $K$ be an object of the derived category $\D^-(\Coh(X\times Y))$.
Let $\pr_1,\pr_2: X\times Y \to X, Y$ be the projections.
For every $\CF \in \D^-(\Coh(X))$ we consider the object
$$
\Phi_K(\CF) = R\pr_{2*}(\pr_1^*\CF \lotimes K) \in \D^-(\Coh(Y)).
$$
The object $\Phi_K(\CF)$ will be called the {\em evaluation} of $K$ on $\CF$.
It is clear that evaluation is functorial, both in $\CF$ and in $K$.
Functoriality in $\CF$ means that $\Phi_K$ is a functor
from the derived category $\D^-(\Coh(X))$ to the derived category $\D^-(\Coh(Y))$
(in other terminology such functors are referred to as {\em integral} or
{\em kernel} functors and the objects $K$ are referred to as {\em kernels}).
Functoriality in $K$ means that to every morphism of kernels $\phi:K_1 \to K_2$
in $\D^-(\Coh(X\times Y))$ corresponds a morphism $\Phi_{K_1}(\CF) \to \Phi_{K_2}(\CF)$
in $\D^-(\Coh(Y))$ which we call {\em evaluation} of $\phi$ on~$\CF$.

\subsection{Atiyah classes}
The Atiyah class was introduced in \cite{At} for the case of vector bundles and
in~\cite{Ill} for any complex of coherent sheaves $\FFF$.
Let $Y$ be an algebraic variety. Let $\Delta:Y \to Y\times Y$ denote the diagonal embedding.
Let $\Delta(Y)^{(2)} \subset Y\times Y$ denote the second infinitesimal neighborhood of the diagonal
$\Delta(Y) \subset Y\times Y$. In other words, if $\CI_\Delta$ is the sheaf of ideals
of the diagonal $\Delta(Y) \subset Y\times Y$, then $\Delta(Y)^{(2)}$ is the closed subscheme
of $Y\times Y$ defined by the sheaf of ideals $\CI_\Delta^2$. Note that
$\CI_\Delta/\CI_\Delta^2 \cong \CN^\dual_{\Delta(Y)/Y\times Y} \cong \Omega_Y$, hence we have
the following exact sequence
\begin{equation}\label{uat}
0 \to \Delta_*\Omega_Y \to \CO_{\Delta(Y)^{(2)}} \to \Delta_*\CO_Y \to 0.
\end{equation}
The corresponding class $\widetilde\At \in \Ext^1(\Delta_*\CO_Y,\Delta_*\Omega_Y)$ is called
the {\em  universal Atiyah class} of $Y$.

Evaluation produces from the universal Atiyah class the usual Atiyah classes of sheaves on $Y$.
Indeed, it is clear that $\Phi_{\Delta_*\CO_Y}(\CF) \cong \CF$,
$\Phi_{\Delta_*\Omega_Y}(\CF) \cong \CF\otimes\Omega_Y$,
so evaluation of $\widetilde{\At}$ on $\CF$ gives a class $\At_\CF\in\Ext^1(\CF,\CF\otimes\Omega_Y)$
which is just the Atiyah class of $\CF$. See \cite{HL}, 10.1.5 for a
representation of $\At_\CF$ by a \v Cech cocycle, or \cite{ALJ}
for an approach via simplicial spaces.

\section{Closed 2-forms on moduli spaces of sheaves}
\label{c2f}

Let $Y$ be a smooth complex projective variety of dimension $n$.
Consider a moduli space $\FM$  of stable sheaves on $Y$ as defined
in Sect. 1 of \cite{S} (see also \cite{LP})
or, more generally, any pre-scheme representing an open part of
the moduli functor ${\bf Spl}_Y$ of simple sheaves \cite{AK}. All the
considerations of this section have also their analytic counterpart
in the case when $Y$ is a compact K\"ahler manifold
and $\FM$ is the analytic moduli space (maybe, non-Hausdorff) of simple
vector bundles on $Y$, see \cite{K}. However, we will need for later
applications the case when $\FM$ is a component of the moduli space
of torsion sheaves, so we cannot restrict ourselves to vector bundles.

For any sheaf $\FFF$
whose isomorphism class $[\FFF ]\in\FM$, the vector
space $\Ext^1 (\FFF ,\FFF)$ is naturally identified
with the tangent space $T_{[\FFF ]}\FM$, and
$\Ext^2 (\FFF ,\FFF)$ is the obstruction space (see
\cite{Muk-1}).
The Yoneda pairing
$$\Ext^i(\FFF ,\FFF )\times\Ext^j(\FFF ,\FFF )
\xrightarrow{\mbox{\scriptsize Yoneda}}
\Ext^{i+j}(\FFF ,\FFF ),\ (a,b)\mapsto a\circ b$$
for $i=j=1$ provides the bilinear map
\begin{equation}\label{yo}
\Lambda :\Ext^1 (\FFF ,\FFF)\times \Ext^1 (\FFF ,\FFF)
\lra \Ext^{2}(\FFF ,\FFF ).
\end{equation}
The obstruction map $\ob :\Ext^1(\FFF ,\FFF )\lra \Ext^{2}(\FFF ,\FFF )$
is expressed in terms of the Yoneda pairing by
$a\mapsto a\circ a$. It has the following sense:
an element $a\in \Ext^1(\FFF ,\FFF )$ defines the isomorphism
class of a flat deformation $\FFF^{(1)}$ of $\FFF$ over $Y\times\Spec
\CC [t]/(t^2)$, and $\ob (a)=0$ if and only if $\FFF^{(1)}$ can be extended
further to a flat deformation $\FFF^{(2)}$ over $Y\times\Spec
\CC [t]/(t^3)$. In particular, if $[\FFF ]$ is a smooth point
of $\FM$, then all the obstructions vanish: $\ob (a)=0\
\forall\ a\in \Ext^1(\FFF ,\FFF )$. Thus we have the following
statement:

\begin{lemma}\label{skew}
Let $\FM^\sm $ denote the smooth locus of $\FM$, and
let $[\FFF ]$ be a point
of $\FM^\sm $. Then the bilinear map $\Lambda$ defined by~\eqref{yo}
is skew symmetric: $\Lambda (a,b)=-\Lambda (b,a)$ for all
$a,b\in \Ext^1(\FFF ,\FFF )$.
\end{lemma}

Buchweitz and Flenner \cite{BuF1}, \cite{BuF2} define the map
\begin{equation}\label{sigma-r}
\sigma=\sum_{q\geq 0}\si_q: \Ext^{2}(\FFF ,\FFF )\lra\bigoplus_{q\geq 0}
H^{q+2}(X,\Omega_X^q),\ \si:c\mapsto \Tr\, (\exp (-\At (\FFF))\circ c),
\end{equation}
which coincides with Bloch's
semiregularity map in the case when $\FFF =\OOO_Z$
for a subscheme $Z$. It is proved in loc. cit. that
for any simple coherent sheaf $\FFF$, the map $\si$
plays the same role for the moduli space of sheaves as
Bloch's semiregularity map for subschemes:
$\FM$ is smooth in~$[\FFF ]$ if $\si$ is injective.
This is the reason to call it semiregularity map for sheaves.

Composing $\Lambda$ with
$
\sigma :\Ext^{2}(\FFF ,\FFF )\lra\bigoplus
H^{q+2}(Y,\Omega_Y^q)$ for $\FFF\in\FM^\sm $,
we obtain a family of skew-symmetric bilinear forms
on $\Ext^{1}(\FFF ,\FFF )=T_{[\FFF ]}\FM$, each one of
which corresponds to some element of the dual space $(\bigoplus
H^{q+2}(Y,\Omega_Y^q))^\dual=
\bigoplus H^{n-q-2}(Y,\Omega_Y^{n-q})$. These forms
fit into exterior 2-forms on $\FM^\sm$ and are linear combinations
of the forms $\alpha_{\omega}$ with $\omega\in H^{n-q-2}(Y,\Omega_Y^{n-q})$,
defined as follows:
\begin{equation}\label{alpha}
\alpha_{\omega}(v_1,v_2) = \Tr(\At_\CF^q \circ v_1 \circ v_2) \ccup \omega \ccap [Y],
\end{equation}
where $v_1,v_2 \in T_{[\FFF ]}\FM$ and $[Y]$ is the fundamental class of $Y$.

\begin{theorem}\label{closed}
Let $Y$ be a smooth complex projective variety of dimension $n$, $\FM$
a moduli space  of stable or simple sheaves on $Y$,
$\omega\in H^{n-q-2}(Y,\Omega_Y^{n-q})$. Then formula
\eqref{alpha} defines a closed  $2$-form $\alpha_{\omega}
\in H^0(\FM^\sm ,\Omega^2)$.
\end{theorem}

\begin{proof}
This statement generalizes the known result for the case when
$Y$ is a surface
(then $\sigma$ is just the trace map with values in $H^2(Y,\OOO_Y)$).
Different proofs in the surface case for moduli of vector bundles
can be found in \cite{Muk-2}, \cite{O'G}, \cite{Bot-1}. A proof
for moduli of sheaves on a surface was given in
\cite{HL}. Our proof, given below, is obtained by a slight modification of the argument from \cite{HL}.

It suffices to prove that given a smooth affine variety $S$, for any
$S$-flat sheaf $\FFF$ on $S\times Y$
defining a classifying morphism $\tau :S\lra\FM^\sm $,
$s\mapsto [\FFF_s]$,
the pullback $\tau^*(\alpha_{\omega})\in
H^0(S ,\Omega^2_S)$
is closed. The latter 2-form is
the following map:
\begin{multline}\label{comp}
T_sS\times T_sS\xrightarrow{KS\times KS}
\Ext^{1}(\FFF_s ,\FFF_s )\times \Ext^{1}(\FFF_s ,\FFF_s )
\xrightarrow{\mbox{\scriptsize Yoneda}}\Ext^{2}(\FFF_s ,\FFF_s )\\
\xrightarrow{\At (\FFF_s )^{q}}
\Ext^{q+2}(\FFF_s ,\FFF_s \otimes\Omega_Y^q)\xrightarrow{\Tr}
H^{q+2}(Y,\Omega_Y^q)\xrightarrow{\cccup\omega}H^{n}(Y,\Omega_Y^n)\simeq\CC\ ,
\end{multline}
where $KS$ stands for the Kodaira--Spencer map.

The Kodaira--Spencer map has the following description in terms of
the Atiyah class $\At_{S \times Y}(\FFF )$. Denote by $A(\FFF)$ the
image of $\At_{S \times Y}(\FFF )$ in $H^0(S ,\EXT^1_{\pr_1}(\FFF
,\FFF\otimes\Omega^1_{S \times Y}))$, where $\pr_i$ is the
projection of $S \times Y$ to the $i$-th factor ($i=1,2$), and
$\EXT^i_{\pr_1}$ stands for the $i$-th derived functor of
$\pr_{1*}\circ\HOM : \D^b(\Coh(S\times Y))^\circ\times\D^b(\Coh(S\times Y)) \to \D^b(\Coh(S))$.
Note that we have a relative analog of the local-to-global spectral sequence
$H^p(S,\EXT^q_{\pr_1}(\CF,\CG)) \impl \Ext^{p+q}(\CF,\CG)$.
Since $S$ is affine this spectral sequence degenerates in the second term,
so that we have an isomorphism $\Ext^i(\CF,\CG) \cong H^0(S,\EXT^i_{\pr_1}(\CF,\CG))$.

Write
$A(\FFF )=A'(\FFF )+A''(\FFF )$ according to the direct sum
decomposition $\Omega^1_{S \times Y} = \pr_1^*\Omega^1_S \oplus
\pr_2^*\Omega^1_Y$. Then $KS$ is the composition
\begin{multline}\label{ks}
\TTT S \xrightarrow{1\otimes A '(\FFF )}\TTT S\otimes
\EXT^1_{\pr_1}(\FFF ,\FFF\otimes\pr_1^*\Omega^1_S )\lra\\
\lra \EXT^1_{\pr_1}(\FFF ,\FFF\otimes\pr_1^*
(\Omega^{1\dual}_S \otimes\Omega^1_S ))\lra
\EXT^1_{\pr_1}(\FFF ,\FFF ).
\end{multline}

Combining (\ref{comp}) and (\ref{ks}), we see that
$$
\tau^*(\alpha_\omega)=(\omega\ccup\gamma)\ccap [Y] ,\ \ \
\gamma =\Tr (A''(\FFF )^q
\circ A'(\FFF )^2) \in H^0(S,\Omega^2_S) \otimes H^{q+2}(Y,\Omega^q_Y).
$$
Consider also the class
$\tilde{\gamma}=\Tr (\At_{S \times Y}(\FFF )^{q+2}) \in
H^{q+2}(S\times Y,\Omega^{q+2}_{S\times Y})$.
According to \cite{HL}, Sect.~10.1.6, it is $d_{S\times Y}$-closed,
where $d_{S\times Y}=d_S\otimes 1+1\otimes d_Y$ is the natural differential
on $H^p(S\times Y, \Omega^q_{S\times Y})$ induced by the De Rham differential
on $\Omega_{S\times Y}^\fatdot$. Hence
$$
d_{S\times Y}(\omega\ccup\tilde{\gamma})=
d_{S\times Y}(\omega )\ccup\tilde{\gamma}+(-1)^{n-q}\omega\ccup d_{S\times Y}(
\tilde{\gamma})=0.
$$
Recall that the groups $H^p(S\times Y, \Omega^q_{S\times Y})$
have a K\"unneth decomposition
$$
H^p(S\times Y, \Omega^q_{S\times Y})=\bigoplus_{i,j}
H^i(S,\Omega_S^j)\otimes H^{p-i}(Y,\Omega_Y^{q-j}).
$$
Since $S$ is affine we have $H^i(S,\Omega_S^j)=0$ for $i>0$,
therefore $\omega\ccup\tilde{\gamma} \in H^n(S\times Y,\Omega^{n+2}_{S\times Y})$
is the sum of the K\"unneth components
$f_j\in H^0(S,\Omega_S^j)\otimes H^{n}(Y,\Omega_Y^{n+2-j})$, $j\geq 2$, and we have
$f_2={q+2 \choose 2}\omega\ccup\gamma  \in H^0(S,\Omega_S^2)\otimes H^{n}(Y,\Omega_Y^{n})$.
As $Y$ is projective, $d_Y$ vanishes on $H^{p-i}(Y,\Omega_Y^{q-j})$,
hence the closedness of $\omega\ccup\tilde{\gamma}$ implies that
of $f_j$ for any $j$.
In particular, $\omega\ccup\gamma$ is closed.
Let us represent $\omega\ccup\gamma$ in the form
$\tau^*(\alpha_\omega)\otimes \eta$ with $\tau^*(\alpha_\omega) \in H^0(S,\Omega_S^2)$, where
$\eta$ is a generator of $H^n(Y,\Omega^n_Y)$, dual to $[Y]$. Then
$$
0=d_{S\times Y}(\omega\ccup\gamma )=d_{S\times Y}(\tau^*(\alpha_\omega)
\otimes \eta )= d_S(\tau^*(\alpha_\omega))\otimes\eta,
$$
which implies that $d_S(\tau^*(\alpha_\omega))=0$.
\end{proof}

Thus we have constructed closed 2-forms $\alpha_{\omega}$
on $\FM^\sm$. In general, these forms may be
degenerate, but, as we will see, they are symplectic in some examples.

\section{The linkage class}

Let $M$ be an algebraic variety and $Y \subset M$, a locally complete intersection
subvariety of codimension $m$. Denote by $i:Y \to M$ the embedding.
Let $\CF$ be a coherent sheaf on $Y$, then $i_*\CF$ is a coherent sheaf
on $M$ supported on $Y$. As we have shown in Lemma~\ref{pbpf},
the derived pullback $Li^*(i_*\CF)$ considered as an object
of the derived category $\D^b(\Coh(Y))$ is a complex with $(m+1)$
nontrivial cohomology, $\CF$ at degree $0$, $\CF\otimes\CN^\dual_{Y/M}$ at degree $-1$,
and so on. Consider the canonical filtration of this object.
Its associated graded factors are the shifted cohomology sheaves,
explicitly $\CF\otimes\wedge^k\CN^\dual_{Y/M}[k]$.
So the object $Li^*i_*\CF$ provides us with extension classes
$\epsilon^k_\CF \in \Ext^1(\CF\otimes\wedge^k\CN^\dual_{Y/M}[k],\CF\otimes\wedge^{k+1}\CN^\dual_{Y/M}[k+1]) \cong
\Ext^2(\CF\otimes\wedge^k\CN^\dual_{Y/M},\CF\otimes\wedge^{k+1}\CN^\dual_{Y/M})$.
The most important of them, $\epsilon_\CF := \epsilon^0_\CF \in \Ext^2(\CF,\CF\otimes\CN^\dual_{Y/M})$
will be called {\em the linkage class} of $\CF$.

If $Y \subset M$ is a divisor, the linkage class
$\epsilon_\CF \in \Ext^2(\CF,\CF\otimes\CN^\dual_{Y/M}) \cong \Ext^2(\CF,\CF\otimes\CO_Y(-Y))$
completely determines the derived pullback $Li^*i_*\CF$,
namely there is a distinguished triangle
\begin{equation}\label{dlc}
\xymatrix@1{Li^*i_*\CF \ar[r] & \CF \ar[r]^-{\epsilon_\CF} & \CF\otimes\CO_Y(-Y)[2]}.
\end{equation}
In other words, $Li^*i_*\CF$, up to a shift, is a cone of $\epsilon_\CF$.

\begin{proposition}
The linkage class $\epsilon_\CF \in \Ext^2(\CF,\CF\otimes\CN^\dual_{Y/M})$
is defined for any object $\CF$ of the derived category $\D^b(\Coh(Y))$
and is functorial, i. e. for any morphism $\varphi:\CF \to \CG$ in $\D^b(\Coh(Y))$
we have a commutative diagram
$$
\xymatrix@C=2cm{
\CF \ar[r]^-{\epsilon_\CF} \ar[d]_\varphi & \CF\otimes\CN^\dual_{Y/M}[2] \ar[d]^{\varphi\otimes 1} \\
\CG \ar[r]^-{\epsilon_\CG} & \CG\otimes\CN^\dual_{Y/M}[2]
}
$$
\end{proposition}
\begin{proof}
Let $\Delta:Y \to Y\times Y$ be the diagonal embedding,
and denote $\tilde\imath = (1\times i):Y\times Y \to Y\times M$.
Then by Lemma~\ref{pbpf}, for any coherent sheaf
$F$ on $Y\times Y$ the cohomology sheaves of the derived pullback
$L\tilde\imath^*\tilde\imath_*F$ are isomorphic to $F\otimes \pr_2^*\wedge^k\CN^\dual_{Y/M}$.
Therefore, we have an extension class
$\tilde\epsilon_F \in \Ext^2(F,F\otimes\pr_2^*\CN^\dual_{Y/M})$.
Take $F = \Delta_*\CO_Y$ and consider $\tilde\epsilon := \tilde\epsilon_{\Delta_*\CO_Y}$
as a morphism $\Delta_*\CO_Y \to \Delta_*\CO_Y\otimes\pr_2^*\CN^\dual_{Y/M}[2] \cong \Delta_*\CN^\dual_{Y/M}[2]$.
Now for any object $\CF \in \D^b(\Coh(Y))$ the evaluation of $\tilde\epsilon_{\Delta_*\CO_Y}$ gives a morphism
\begin{equation}\label{dereps}
\xymatrix@1{
\CF \cong \pr_{2*}(\pr_1^*\CF \otimes \Delta_*\CO_Y) \ar[r]^-{\tilde\epsilon} &
\pr_{2*}(\pr_1^*\CF \otimes \Delta_*\CN^\dual_{Y/M}[2]) \cong \CF\otimes\CN^\dual_{Y/M}[2].
}
\end{equation}
It is clear that if $\CF$ is a coherent sheaf, then this morphism considered
as an element of $\Ext^2(\CF,\CF\otimes\CN^\dual_{Y/M})$ coincides with
the linkage class $\epsilon_\CF$ defined above, so~\eqref{dereps}
can be considered as an extension of the definition of the linkage class
to the whole derived category. Moreover, \eqref{dereps} also shows that $\epsilon_\CF$ is
the evaluation of the universal class $\tilde\epsilon$ and hence is functorial.
\end{proof}

Actually, the linkage class can be expressed in terms of the Atiyah class.
Consider the adjunction exact sequence
\begin{equation}\label{adjseq}
\xymatrix@1{0 \ar[r]& \CN^\dual_{Y/M} \ar[r]^{\kappa}& \Omega_{M|Y} \ar[r]^{\rho}& \Omega_Y \ar[r]& 0}
\end{equation}
and denote by $\kappa = \kappa_{Y/M}:\CN^\dual_{Y/M} \to \Omega_{M|Y}$,
$\rho = \rho_{Y/M}:\Omega_{M|Y} \to \Omega_Y$ the maps in~\eqref{adjseq},
and by $\nu = \nu_{Y/M} \in \Ext^1(\Omega_Y,\CN^\dual_{Y/M})$
the extension class of~\eqref{adjseq}.

\begin{theorem}\label{epsat}
Let $i:Y \to M$ be a locally complete intersection.

$(i)$ For any $\CF \in \D^b(\Coh(Y))$,
the linkage class $\epsilon_\CF \in \Ext^2(\CF,\CF\otimes\CN^\dual_{Y/M})$
is the product of the Atiyah class $\At_\CF \in \Ext^1(\CF,\CF\otimes\Omega_Y)$
with $\nu_{Y/M}$. In other words $\epsilon_\CF = (1_\CF\otimes\nu_{Y/M})\circ\At_\CF$.

$(ii)$ For any $\CG \in D^b(\Coh(M))$ we have
$$
\At_{Li^*\CG} \cong \rho_*((\At_\CG)_{|Y}),
$$
where $\rho_*:\Ext^1(Li^*\CG,Li^*\CG\otimes\Omega_{M|Y}) \to \Ext^1(Li^*\CG,Li^*\CG\otimes\Omega_Y)$
is the pushout via $\rho:\Omega_{M|Y} \to \Omega_Y$.

$(iii)$ For any $\CF \in D^b(\Coh(Y))$ the image of the Atiyah class
$\At_{i_*\CF} \in \Ext^1(i_*\CF,i_*\CF\otimes\Omega_M)$ in
$\Hom (L_1 i^*i_*\CF,\CF\otimes\Omega_{M|Y})=
H^0(M,i_*(\CF^\dual\otimes\CF\otimes\CN_{Y/M}\otimes\Omega_{M|Y})) =
\Hom(\CF\otimes\CN_{Y/M}^\dual,\CF\otimes\Omega_{M|Y})$
%
%
equals $1_\CF\otimes\kappa$. In the particular case when $\CF$ is locally
free on $Y$, we have $\HOM(\CF\otimes\CN_{Y/M}^\dual,\CF\otimes\Omega_{M|Y})=
\EXT^1(i_*\CF,i_*\CF\otimes\Omega_M)$, and the image of $\At_{i_*\CF}$
in $H^0(M,\EXT^1(i_*\CF,i_*\CF\otimes\Omega_M))$ under tha map coming from
the local-to-global spectral sequence is also $1_\CF\otimes\kappa$.
\end{theorem}
\begin{proof}
Let $\Delta_Y:Y \to Y\times Y$ and $\Delta_M:M \to M\times M$ denote the diagonal
embeddings, and let $\Gamma:Y \to M\times Y$ be the graph of $i:Y \to M$.
Then we have a commutative diagram
$$
\xymatrix{
Y \ar@{=}[r] \ar[d]_{\Delta_Y} &
Y \ar[r]^i \ar[d]_{\Gamma} &
M \ar[d]^{\Delta_M} \\
Y\times Y \ar[r]^{i\times 1} &
M\times Y \ar[r]^{1\times i} &
M\times M
}
$$
Note that both squares of the diagram are cartesian.
Let $\Delta_M(M)^{(2)}$, $\Gamma(Y)^{(2)}$ and $\Delta_Y(Y)^{(2)}$
denote the second infinitesimal neighborhoods of $\Delta_M(M)$ in $M\times M$,
$\Gamma(Y)$ in $M\times Y$ and $\Delta_Y(Y)$ in $Y\times Y$ respectively.
On $M\times M$ we have the short exact sequence
$$
0 \to \Delta_{M*}\Omega_M \to \CO_{\Delta_M(M)^{(2)}} \to \Delta_{M*}\CO_M \to 0
$$
representing the universal Atiyah class on $M$.
Consider its pullback to $M\times Y$. Since $M\times Y$ intersects $\Delta_M(M)$
transversely, the higher derived inverse images of $1\times i$ are zero, and we have an isomorphism
$(1\times i)^*\Delta_{M*}\CG \cong \Gamma_*i^*\CG$ for any coherent sheaf $\CG$ on $M$.
Moreover, it is clear that $(1\times i)^*(\CO_{\Delta_M(M)^{(2)}}) \cong \CO_{\Gamma(Y)^{(2)}}$,
hence we obtain the exact sequence
\begin{equation}\label{pbmy}
0 \to \Gamma_*\Omega_{M|Y} \to \CO_{\Gamma(Y)^{(2)}} \to \Gamma_*\CO_Y \to 0.
\end{equation}
By definition, this sequence represents the restriction to $Y$
of the universal Atiyah class of $M$. On the other hand, it is clear
that the map $\CO_{\Gamma(Y)^{(2)}} \to \Gamma_*\CO_Y$ factors through
$(i\times 1)_*\CO_{\Delta_Y(Y)^{(2)}}$, hence we have a commutative diagram
$$
\xymatrix@R=10pt{
0 \ar[r] & \Gamma_*\Omega_{M|Y} \ar[r]\ar[d]_\rho& \CO_{\Gamma(Y)^{(2)}} \ar[r]\ar[d]& \Gamma_*\CO_Y \ar[r]\ar@{=}[d]& 0\\
0 \ar[r] & (i\times 1)_*\Delta_*\Omega_Y \ar[r]& (i\times 1)_*\CO_{\Delta_Y(Y)^{(2)}} \ar[r]& (i\times 1)_*\Delta_*\CO_Y \ar[r]& 0
}
$$
Since the bottom line represents the universal Atiyah class of $Y$, we deduce part~$(ii)$.

Now consider the pullback of~\eqref{pbmy} to $Y\times Y$.
This time the intersection of $Y\times Y$ with $\Gamma(Y)$ is not transversal.
Actually, we have $\Gamma_*\CF \cong (1\times i)_*\Delta_{Y_*}\CF$
for any coherent sheaf $\CF$ on~$Y$, hence by Lemma~\ref{pbpf} we have
$(1\times i)^*\Gamma_*\CF \cong \Delta_{Y*}\CF$,
$L_1(1\times i)^*\Gamma_*\CF \cong \Delta_{Y*}(\CF\otimes\CN^\dual_{Y/M})$, and so on.
Moreover, it is clear that $(i\times 1)^*(\CO_{\Gamma(Y)^{(2)}}) \cong \CO_{\Delta_Y(Y)^{(2)}}$,
hence we obtain the following long exact sequence
\begin{equation}\label{pbyy}
\dots \to
L_1(1\times i)^*\CO_{\Gamma(Y)^{(2)}} \to \Delta_{Y*}\CN^\dual_{Y/M} \to
\Delta_{Y*}\Omega_{M|Y} \to \CO_{\Delta_Y(Y)^{(2)}} \to \Delta_{Y*}\CO_Y \to 0.
\end{equation}
The map $\Delta_{Y*}\CN^\dual_{Y/M} \to \Delta_{Y*}\Omega_{M|Y}$
in this sequence is the image
of the universal Atiyah class of $M$, restricted to $Y$,
under the natural map
\begin{multline*}
\Ext^1(\Gamma_*\CO_Y,\Gamma_*\Omega_{M|Y}) \cong
\Ext^1((i\times 1)_*\Delta_*\CO_Y,(i\times 1)_*\Delta_*\Omega_{M|Y}) \cong \\ \cong
\Ext^1(L(i\times 1)^*(i\times 1)_*\Delta_*\CO_Y,\Delta_*\Omega_{M|Y}) \to
\Hom(L_1(i\times 1)^*(i\times 1)_*\Delta_*\CO_Y,\Delta_*\Omega_{M|Y}) \to \\ \to
\Hom((\CN^\dual_{Y/M}\boxtimes\CO_Y)\otimes\Delta_*\CO_Y,\Delta_*\Omega_{M|Y}) \cong
\Hom(\Delta_*\CN^\dual_{Y/M},\Delta_*\Omega_{M|Y}),
\end{multline*}
so for part~$(iii)$ it suffices to check that this map is $\kappa$.
Comparing~\eqref{pbyy} with the sequence
$$
0 \to \Delta_{Y*}\Omega_Y \to \CO_{\Delta_Y(Y)^{(2)}} \to \Delta_{Y*}\CO_Y \to 0
,$$
we see that the map $\Delta_{Y*}\Omega_{M|Y} \to \CO_{\Delta_Y(Y)^{(2)}}$
in~\eqref{pbyy} factors through $\Delta_{Y*}\Omega_Y$. It is clear that
the arising map $\Omega_{M|Y} \to \Omega_Y$ is the restriction of differential forms,
hence its kernel is isomorphic to $\CN^\dual_{Y/M}$. Thus we see that the map
$L_1(1\times i)^*\CO_{\Gamma(Y)^{(2)}} \to \Delta_{Y*}\CN^\dual_{Y/M}$ in~\eqref{pbyy}
must be zero and the last 4 terms of~\eqref{pbyy} form an exact sequence
\begin{equation}\label{pbyy1}
\xymatrix@1{
0 \ar[r]& \Delta_{Y*}\CN^\dual_{Y/M} \ar[r]^-\kappa& \Delta_{Y*}\Omega_{M|Y} \ar[r]&
\CO_{\Delta_Y(Y)^{(2)}} \ar[r]& \Delta_{Y*}\CO_Y \ar[r]& 0.}
\end{equation}
So, part~$(iii)$ follows.
Finally, the Yoneda class of the extension~\eqref{pbyy1} by definition equals
$\tilde\epsilon \in \Ext^2(\Delta_{Y*}\CO_Y,\Delta_{Y_*}\CN^\dual_{Y/M})$.
On the other hand, as we have seen above this class factors as
$\tilde\epsilon = \Delta_{Y*}(\nu_{Y/M}) \circ \widetilde\At$,
where $\widetilde\At \in \Ext^1(\Delta_{Y*}\CO_Y,\Delta_{Y*}\Omega_Y)$
is the universal Atiyah class. Evaluating this equality on any $\CF \in \D^b(\Coh(Y))$
we deduce part $(i)$ of the theorem.
\end{proof}

\section{Application to cubic fourfold}\label{sc4}

From now on we take $M = \PP^5$ and $Y \subset \PP^5$, a smooth cubic fourfold.
Let $i:Y \to \PP^5$ denote the embedding.
Note that $\Omega^4_Y \cong \CO_Y(-3) \cong \CN^\dual_{Y/\PP^5}$.
This coincidence will be important below.
Consider the full triangulated subcategory $\CCC_Y \subset \D^b(\Coh(Y))$
defined by
\begin{equation}\label{ccy}
\begin{array}{rcl}
\CCC_Y &=&
\{\CF \in \D^b(\Coh(Y))\ |\
\Ext^\bullet(\CO_Y,\CF) =
\Ext^\bullet(\CO_Y(1),\CF) =
\Ext^\bullet(\CO_Y(2),\CF) = 0
\} \\
&=&
\{\CF \in \D^b(\Coh(Y))\ |\
H^\bullet(Y,\CF) =
H^\bullet(Y,\CF(-1)) =
H^\bullet(Y,\CF(-2)) = 0
\}.
\end{array}
\end{equation}

\begin{proposition}\label{lc-on-exts}
Assume that $\CF,\CG \in \CCC_Y$.
Then the multiplication by the linkage class
$\epsilon_\CG \in \Ext^2(\CG,\CG(-3))$ induces an isomorphism
$\Ext^p(\CF,\CG) \cong \Ext^{p+2}(\CF,\CG(-3))$ for all $p$.
\end{proposition}
\begin{proof}
Consider the Beilinson spectral sequence for $i_*\CG$ (see \cite{Bei,OSS})
$$
E_1^{-p,q} = H^q(\PP^5,i_*\CG(-p))\otimes \Omega_{\PP^5}^p(p) \Longrightarrow i_*\CG
\qquad (p = 0,1,\dots,5).
$$
Note that $H^\bullet(\PP^5,i_*\CG(-p)) = H^\bullet(Y,\CG(-p))$,
hence $E_1^{0,q} = E_1^{-1,q} = E_1^{-2,q} = 0$ for all $q$ since $\CG \in \CCC_Y$.
It follows that the derived pullback $Li^*i_*\CG$ is contained in the triangulated
subcategory of $\D^b(\Coh(Y))$ generated by $i^*\Omega_{\PP^5}^3(3)$,
$i^*\Omega_{\PP^5}^4(4)$, and $i^*\Omega_{\PP^5}^5(5)$. On the other hand,
the standard resolutions
$$
\begin{array}{r}
0 \to \CO_{\PP^5}(-3) \to \CO_{\PP^5}(-2)^{\oplus 6} \to \CO_{\PP^5}(-1)^{\oplus 15} \to \Omega_{\PP^5}^3(3) \to 0,\\[2pt]
0 \to \CO_{\PP^5}(-2) \to \CO_{\PP^5}(-1)^{\oplus 6} \to \Omega_{\PP^5}^4(4) \to 0,
\end{array}
$$
and an isomorphism $\CO_{\PP^5}(-1) \cong \Omega_{\PP^5}^5(5)$ show that
this subcategory coincides with the subcategory of $\D^b(\Coh(Y))$ generated
by $\CO_Y(-3)$, $\CO_Y(-2)$ and $\CO_Y(-1)$. But note that the Serre duality gives
$$
\begin{array}{l}
\Ext^p(\CF,\CO_Y(-3)) \cong \Ext^{4-p}(\CO_Y,\CF)^\dual = 0,\\[2pt]
\Ext^p(\CF,\CO_Y(-2)) \cong \Ext^{4-p}(\CO_Y(1),\CF)^\dual = 0,\\[2pt]
\Ext^p(\CF,\CO_Y(-1)) \cong \Ext^{4-p}(\CO_Y(2),\CF)^\dual = 0,
\end{array}
$$
which implies that $\Ext^\bullet(\CF,Li^*i_*\CG) = 0$.
Applying the functor $\Hom(\CF,-)$ to the distinguished triangle~\eqref{dlc} for
$\CG$, we deduce the proposition.
\end{proof}

\begin{remark}
Combining the isomorphism $\Ext^p(\CF,\CG) \cong \Ext^{p+2}(\CF,\CG(-3))$ of the proposition
with the Serre duality $\Ext^{p+2}(\CF,\CG(-3)) \cong \Ext^{2-p}(\CG,\CF)^\dual$, we obtain
a duality
\begin{equation}\label{sd_ccc}
\Ext^p(\CF,\CG) \cong \Ext^{2-p}(\CG,\CF)^\dual
\qquad \text{for any \ \ $\CF,\CG \in \CCC_Y$}.
\end{equation}
This duality, in fact, is the Serre duality for the triangulated category $\CCC_Y$.
In other words, the Serre functor (see~\cite{BK}) of $\CCC_Y$ equals to the shift by 2 functor.
This fact was proved earlier in~\cite{Ku1} by the same argument.
\end{remark}

Using the isomorphism
$\Ext^1(\Omega_Y,\CN^\dual_{Y/\PP^5}) \cong H^1(Y,\CT_Y\otimes\CO_Y(-3)) \cong H^1(Y,\Omega^3_Y)$,
we consider the extension class $\nu_{Y/\PP^5} \in \Ext^1(\Omega_Y,\CN^\dual_{Y/\PP^5})$
of the adjunction sequence~\eqref{adjseq} as an element of $H^1(Y,\Omega^3_Y)$.

\begin{theorem}\label{sf1}
Let $\FM$ be a moduli space of stable sheaves on a cubic $4$-fold $Y$ such that for every sheaf $\CF$
with $[\CF]\in\FM$ we have $H^\bullet(Y,\CF) = H^\bullet(Y,\CF(-1)) = H^\bullet(Y,\CF(-2)) = 0$.
Then the closed $2$-form $\alpha_\nu \in H^0(\FM^\sm,\Omega^2)$ corresponding to the class
$\nu = \nu_{Y/\PP^5} \in H^1(Y,\Omega^3_Y)$
is nondegenerate.
\end{theorem}
\begin{proof}
Let $[\CF] \in \FM^\sm$.
Recall that for any $a,b \in \Ext^1(\CF,\CF)$ the form $\alpha_\nu(a,b)$ is defined as
$\nu\circ\Tr(\At_\CF\circ a\circ b)$. By multiplicativity of the trace, this is equal to
$\Tr((1_\CF\otimes\nu)\circ \At_\CF\circ a\circ b)$. In other words, we apply
the Yoneda multiplication map
$$
\Ext^1(\CF\otimes\Omega_Y,\CF\otimes\CN^*_{Y/\PP^5})\otimes
\Ext^1(\CF,\CF\otimes\Omega_Y)\otimes
\Ext^1(\CF,\CF)\otimes\Ext^1(\CF,\CF) \to
\Ext^4(\CF,\CF\otimes\CN^*_{Y/\PP^5})
$$
to $(1_\CF\otimes\nu)\otimes \At_\CF\otimes a\otimes b$
and then the trace map
$$
\Tr:\Ext^4(\CF,\CF\otimes\CN^\dual_{Y/\PP^5}) \cong \Ext^4(\CF,\CF(-3)) \to H^4(Y,\CO_Y(-3)) = \CC.
$$
Since the Yoneda multiplication is associative we have
$$
\Tr((1_\CF\otimes\nu)\circ \At_\CF\circ a\circ b) = \Tr((\epsilon_\CF\circ a)\circ b)
$$
by Theorem~\ref{epsat}~(i). It remains to note that the map $\Ext^1(\CF,\CF) \to \Ext^3(\CF,\CF(-3))$,
$a\mapsto \epsilon_\CF\circ a$ is an isomorphism by Proposition
\ref{lc-on-exts}, and that the Serre duality pairing
$\xymatrix@1{\Ext^3(\CF,\CF(-3))\otimes\Ext^1(\CF,\CF) \ar[r]^-\circ & \Ext^4(\CF,\CF(-3)) \ar[r]^-{\Tr} & H^4(Y,\CO_Y(-3)) = \CC}$
is nondegenerate.
\end{proof}

There are two well-known examples of symplectic moduli spaces of sheaves on a cubic fourfold~$Y$.
The first one \cite{BD} is the Hilbert scheme of lines on $Y$, which we will denote by $F(Y)$.
It was shown in \cite{BD} that $F(Y)$ is an irreducible symplectic variety of dimension 4.
The second one \cite{MT2} is (an open subset of) the moduli space of torsion sheaves
of the form $i_*\CE$, where $\CE$ is a vector bundle of rank $2$
with $c_1 = 0$ and $c_2 = 2[\ell]$ on a smooth hyperplane section
$Y'$ of $Y$, and $i:Y' \to Y$ is the embedding.
This variety $\FMY$ is 10-dimensional and the map $\FMY \to \check{\PP}^5$
taking a sheaf $i_*\CE$ to its support hyperplane section $Y' \subset Y$,
considered as a point of the dual projective space $\check{\PP}^5$
was shown in~\cite{MT2} to be a Lagrangian fibration.

We will show that our results provide constructions of symplectic forms on both
these moduli spaces. The 10-dimensional moduli space $\FMY$ can be dealt with
in a straightforward way. We will explain in Section~\ref{sdim10} that
all sheaves $i_*\CE$ belonging to this moduli space are contained
in the subcategory $\CCC_Y$ of $\D^b(\Coh(Y))$, hence Theorem~\ref{sf1}
applies and gives a symplectic form on~$\FMY$. The details of the construction
can be found in Section~\ref{sdim10}.

The case of the variety $F(Y)$ of lines on $Y$ is slightly more complicated.
For technical reasons it is more convenient to consider $F(Y)$ as
the moduli space of twisted ideal sheaves $\CI_\ell(1)$, where $\ell \subset Y$ is a line.
Certainly, the sheaves $\CI_\ell(1)$ are not contained in $\CCC_Y$
(nor any other twist of $\CO_\ell$ or of $\CI_\ell$).
However, it turns out that a simple endofunctor $\LL:\D^b(\Coh(Y)) \to \D^b(\Coh(Y))$
(the left mutation in $\CO_Y$) takes $\CI_\ell(1)$ to the subcategory
$\CCC_Y$ for every line $\ell$. Explicitly, $\LL(\CI_\ell(1))$ is just
the ``second syzygy sheaf'' of $\ell$ twisted by~$1$, that is
the kernel of the natural map $\CO_Y^{\oplus 4} = H^0(Y,\CI_\ell(1))\otimes\CO_Y \to \CI_\ell(1)$.
In other words, $\CF_\ell := \LL(\CI_\ell(1))$ is the reflexive sheaf
on $Y$ defined by the following exact sequence
\begin{equation}\label{def_fl}
0 \to \CF_\ell \to \CO_Y^{\oplus 4} \to \CO_Y(1) \to \CO_\ell(1) \to 0.
\end{equation}
It is easy to see that for each $\ell$ the sheaf $\CF_\ell$ is stable
and is contained in the subcategory $\CCC_Y$ of $\D^b(\Coh(Y))$.
Therefore, Theorem~\ref{sf1} gives a symplectic form on the module space
$F'(Y)$ of stable sheaves containing sheaves $\CF_\ell$.
We will show in Section~\ref{sdim4} that the map \mbox{$\LL:F(Y) \to F'(Y)$},
$[\ell] \mapsto [\CF_\ell]$ is an open embedding (so $F(Y)$, being
projective, is identified with a connected component of $F'(Y)$)
hence the symplectic form on $F'(Y)$ restricts to a symplectic
form on $F(Y)$. Moreover, we will show that this form coincides
with the form $\alpha_\nu$ defined in Section~\ref{c2f}.

\begin{remark}
Another example of a symplectic variety associated to a cubic fourfold $Y$
was constructed recently by Iliev and Manivel \cite{IMan}.
Unfortunately, we do not know
whether it is possible to realize it as a moduli space of sheaves on $Y$.
It would be interesting to find such a realization. Then Theorem~\ref{sf1}
would give a construction of a symplectic form on this moduli space.
\end{remark}

\section{The variety of lines}\label{sdim4}

Let $Y \subset \PP^5$ be a smooth cubic 4-fold.
Let $F(Y)$ be the Hilbert scheme of lines on $Y$.
We consider $F(Y)$ as the moduli space of sheaves $\CI_\ell(1)$
where $\CI_\ell \subset \CO_Y$ is the ideal of a line $\ell$.

Consider the functor $\LL:\D^b(\Coh(Y)) \to \D^b(\Coh(Y))$ defined as follows
$$
\LL(F) = \mathsf{Cone}\{ \xymatrix@1{H^\bullet(Y,F)\otimes \CO_Y \ar[r]^-{\mathsf{ev}} & F} \},
$$
Here $\mathsf{ev}$ stands for the evaluation homomorphism and $\mathsf{Cone}$ stands for the cone
of a morphism in the derived category. The functor $\LL$ actually is the left mutation
through the exceptional line bundle $\CO_Y$, see~\cite{Bon,GR}.

\begin{lemma}
Let $\ell \subset Y$ be a line. Then $\LL(\CI_\ell(1))[-1]$ is isomorphic
in $\D^b(\Coh(Y))$ to a reflexive sheaf $\CF_\ell$ of rank $3$ on $Y$,
which fits into the exact sequence~\eqref{def_fl}. Moreover, $\CF_\ell \in \CCC_Y$.
\end{lemma}
\begin{proof}
We have $H^\bullet(Y,\CI_\ell(1)) = \CC^4$, hence
$\LL(\CI_\ell(1)) = \mathsf{Cone}\{\xymatrix@1{\CO_Y^{\oplus 4} \ar[r]^-{\mathsf{ev}} & \CI_\ell(1)}\}$.
Note that the space of global sections of $\CO_Y^{\oplus 4}$ here is spanned
by linear functions on $\PP^5$ vanishing on $\ell$, and the map $\CO_Y^{\oplus 4} \to \CI_\ell(1)$
is induced by considering these functions as sections of $\CI_\ell(1)$.
Since the sheaf of ideals of a line is generated by these linear functions,
the evaluation homomorphism is surjective, hence
$\LL(\CI_\ell(1)) = \CF_\ell[1]$, where
$\CF_\ell = \mathsf{Ker}\{\xymatrix@1{\CO_Y^{\oplus 4} \ar[r]^-{\mathsf{ev}} & \CI_\ell(1)}\}$.
Combining the sequence
\begin{equation}\label{flil}
0 \to \CF_\ell \to \CO_Y^{\oplus 4} \to \CI_\ell(1) \to 0
\end{equation}
with the sequence $0 \to \CI_\ell \to \CO_Y \to \CO_\ell \to 0$ twisted by $\CO_Y(1)$,
we deduce~\eqref{def_fl}.

Moreover, using~\eqref{def_fl} to compute $H^\bullet(Y,\CF(-q))$ for $q = 0,1,2$,
we conclude that $\CF_\ell \in \CCC_Y$.
\end{proof}

\begin{proposition}
For any line $\ell \subset Y$, the sheaf $\CF_\ell$ is stable.
Moreover, for $\ell \ne \ell'$ we have $\CF_{\ell} \not\cong \CF_{\ell'}$.
\end{proposition}
\begin{proof}
The sheaf $\CF_\ell$ is reflexive of rank 3 with $c_1(\CF_\ell) = -1$, hence for stability
it suffices to check that $H^0(Y,\CF_\ell) = H^0(Y,\CF_\ell^\dual(-1)) = 0$.
But $H^0(Y,\CF_\ell) = 0$ since $\CF_\ell \in \CCC_Y$, and by Serre duality
$H^0(Y,\CF_\ell^\dual(-1)) = H^4(Y,\CF_\ell(-2)) = 0$
since $\CF_\ell \in \CCC_Y$.

Further, note that~\eqref{def_fl} implies that
$\EXT^1(\CF_\ell,\CO_Y) \cong \EXT^3(\CO_\ell(1),\CO_Y) \cong \CO_\ell$,
whereof it follows that $\CF_\ell \not\cong \CF_{\ell'}$ for $\ell \ne \ell'$.
\end{proof}

\begin{corollary}\label{extfl}
For any line $\ell \subset Y$, we have $\dim\Hom(\CF_\ell,\CF_\ell) = \dim \Ext^2(\CF_\ell,\CF_\ell) = 1$,
$\dim \Ext^1(\CF_\ell,\CF_\ell) = 4$.
\end{corollary}
\begin{proof}
The equality $\dim \Hom(\CF_\ell,\CF_\ell) = 1$ follows from the stability of $\CF_\ell$, and
$\dim \Ext^2(\CF_\ell,\CF_\ell) = 1$ follows from~\eqref{sd_ccc}.
It also follows from~\eqref{sd_ccc} that $\Ext^p(\CF_\ell,\CF_\ell) = 0$
for $p > 2$. Therefore, $\dim\Ext^1(\CF_\ell,\CF_\ell)$ can be computed
by Riemann--Roch.
\end{proof}

\begin{proposition}
The map $\LL:\Ext^1(\CI_\ell(1),\CI_\ell(1)) \to \Ext^1(\CF_\ell,\CF_\ell)$
induced by the functor $\LL$ is an isomorphism.
\end{proposition}
\begin{proof}
Applying the functor $\Hom(-,\CF_\ell)$ to~\eqref{flil} and taking into account
that $\Hom(\CO_Y,\CF_\ell) = 0$, we conclude that
$$
\Ext^p(\CF_\ell,\CF_\ell) \cong \Ext^{p+1}(\CI_\ell(1),\CF_\ell)
$$
for all $p$. Further, note that by Serre duality we have
$$
\Ext^p(\CI_\ell(1),\CO_Y) \cong
\Ext^{4-p}(\CO_Y,\CI_\ell(-2))^\dual \cong
H^{3-p}(\ell,\CO_\ell(-2))^\dual =
\begin{cases}\CC, & \text{if $p = 2$}\\0, & \text{otherwise}\end{cases}
$$
Therefore, applying the functor $\Hom(\CI_\ell(1),-)$ to~\eqref{flil},
we obtain the exact sequence
$$
\arraycolsep=1pt
\begin{array}{rrrrcrrrr}
0 &\to& \Hom(\CI_\ell(1),\CF_\ell) &\to& 0 &\to& \Hom(\CI_\ell(1),\CI_\ell(1)) &\to& \\
&\to& \Ext^1(\CI_\ell(1),\CF_\ell) &\to& 0 &\to& \Ext^1(\CI_\ell(1),\CI_\ell(1)) &\to& \\
&\to& \Ext^2(\CI_\ell(1),\CF_\ell) &\to& \CC^4 &\to& \Ext^2(\CI_\ell(1),\CI_\ell(1)) &\to& \\
&\to& \Ext^3(\CI_\ell(1),\CF_\ell) &\to& 0 &\to& \Ext^3(\CI_\ell(1),\CI_\ell(1)) &\to& \\
&\to& \Ext^4(\CI_\ell(1),\CF_\ell) &\to& 0 &\to& \Ext^4(\CI_\ell(1),\CI_\ell(1)) &\to& 0.
\end{array}
$$
The composition of the map $\Ext^1(\CI_\ell(1),\CI_\ell(1)) \to \Ext^2(\CI_\ell(1),\CF_\ell)$
with the isomorphism $\Ext^2(\CI_\ell(1),\CF_\ell) \cong \Ext^1(\CF_\ell,\CF_\ell)$
clearly coincides with the map induced by the functor $\LL$, hence
$\LL:\Ext^1(\CI_\ell(1),\CI_\ell(1)) \to \Ext^1(\CF_\ell,\CF_\ell)$
is injective. But $\dim\Ext^1(\CI_\ell(1),\CI_\ell(1)) = 4$, since
this is the tangent space to the smooth 4-dimensional moduli space $F(Y)$,
and $\dim\Ext^1(\CF_\ell,\CF_\ell) = 4$ by Corollary~\ref{extfl}. Hence the map $\LL$
is an isomorphism.
\end{proof}

Let $F'(Y)$ denote the moduli space of stable sheaves on $Y$ containing sheaves $\CF_\ell$.
Consider the map $F(Y) \to F'(Y)$ defined by the functor $\LL$, $\ell \mapsto \LL(\CI_\ell(1))[-1] = \CF_\ell$.

\begin{proposition}
The map $\LL:F(Y) \to F'(Y)$, $\ell \mapsto \CF_\ell$ is an isomorphism of $F(Y)$
with a connected component of $F'(Y)$.
\end{proposition}
\begin{proof}
We already know that $\LL$ induces an isomorphism on tangent spaces, hence it is \'etale.
On the other hand, if $\ell \ne \ell'$ then $\CF_\ell \not\cong \CF_{\ell'}$.
Hence $\LL$ is injective. Thus $\LL$ has to be an open embedding.
Since $F(Y)$ is a projective variety, its image is closed.
Therefore, $\LL$ is an isomorphism onto a connected component.
\end{proof}

It is an interesting question, whether $F'(Y) = F(Y)$, or not.

\begin{theorem}
The map $\LL:F(Y) \to F'(Y)$, $\ell \mapsto \CF_\ell$ agrees up to a sign with the forms $\alpha_\nu$.
In particular, the form $\alpha_\nu$ on $F(Y)$ is symplectic.
\end{theorem}
\begin{proof}
Let us denote the form $\alpha_\nu$ on $F(Y)$ by $\alpha$,
and the form $\alpha_\nu$ on $F'(Y)$ by $\alpha'$.
Take any line $\ell$ on $Y$ and any $a,b \in \Ext^1(\CI_\ell(1),\CI_\ell(1))$.
As it was shown in the proof of Theorem~\ref{sf1} the value $\alpha(a,b)$
is the trace of the following composition of morphisms
$\xymatrix@1{\CI_\ell(1) \ar[r]^-b & \CI_\ell(1)[1] \ar[r]^-a & \CI_\ell(1)[2] \ar[r]^-{\epsilon_{\CI_\ell(1)}} & \CI_\ell(-2)[4]}$
in the derived category $\D^b(\Coh(Y))$, and $\alpha'(\LL(a),\LL(b))$
is the trace of the composition
$\xymatrix@1{\CF_\ell \ar[r]^-{\LL(b)} & \CF_\ell[1] \ar[r]^-{\LL(a)} & \CF_\ell[2] \ar[r]^-{\epsilon_{\CF_\ell}} & \CF_\ell(-3)[4]}$.
By functoriality of the Yoneda multiplication and of the linkage class we have the following
commutative diagram in $\D^b(\Coh(Y))$
$$
\xymatrix@R=15pt{
\CF_\ell \ar[r] \ar[d]^{\LL(b)} & \CO_Y^{\oplus 4} \ar[r] \ar[d]^0 & \CI_\ell(1) \ar[r] \ar[d]^b & \CF_\ell[1] \ar[d]^{-\LL(b)}\\
\CF_\ell[1] \ar[r] \ar[d]^{\LL(a)} & \CO_Y^{\oplus 4}[1] \ar[r] \ar[d]^0 & \CI_\ell(1)[1] \ar[r] \ar[d]^a & \CF_\ell[2] \ar[d]^{-\LL(a)} \\
\CF_\ell[2] \ar[r] \ar[d]^{\epsilon_{\CF_\ell}} &
\CO_Y^{\oplus 4}[2] \ar[r] \ar[d]^0 & \CI_\ell(1)[2] \ar[r] \ar[d]^{\epsilon_{\CI_\ell(1)}} & \CF_\ell[3] \ar[d]^{-\epsilon_{\CF_\ell}} \\
\CF_\ell(-3)[4] \ar[r] & \CO_Y(-3)^{\oplus 4}[4] \ar[r] & \CI_\ell(-2)[1] \ar[r] & \CF_\ell(-3)[5]
}
$$
Since the trace is additive and the trace of the second column is 0,
we conclude that $\alpha'(\LL(a),\LL(b)) = -\alpha(a,b)$.
\end{proof}

\begin{remark}
The same argument shows that the 2-form on $F(Y)$ considered as the moduli space
of ideal sheaves $\CI_\ell$ agrees up to a sign with the 2-form on $F(Y)$
considered as the moduli space of structure sheaves $\CO_\ell$.
\end{remark}

In the next section we will compute the form $\alpha$ on $F(Y)$ explicitly in coordinates.

\section{Closed forms on Hilbert schemes}

Let $Y$ be a projective variety. Fix an ample line bundle
$\OOO(1)$ and a Hilbert polynomial $h$ of some reduced
equidimensional proper closed subscheme $Z_0\subset Y$.
Consider the moduli space $\CP$ of stable sheaves on $Y$ with Hilbert polynomial $h$. Then it has an open subscheme $\PPP_0$ parameterizing
torsion free sheaves of rank $1$ on subschemes $Z \subset Y$
with the same Hilbert polynomial as $Z_0$.
We obtain a morphism $p:\CP_0 \to \Hilb^h(Y)$ to the Hilbert scheme
parameterizing the subschemes $Z \subset Y$ with Hilbert polynomial
$h$. The fiber of $p:\CP_0 \to \Hilb^h(Y)$
over a point $[Z]\in\Hilb^h(Y)$, if nonempty, is a partial compactification of
the generalized Picard scheme $\Pic^0(Z)$.

In this section we will give another interpretation of the forms $\alpha_\omega$
constructed in~\eqref{alpha} in the special case when the moduli space is $\CP_0$.
Actually we will show that for some $\omega$ the form $\alpha_\omega$ is the pullback
of a $2$-form on $\Hilb^h(Y)$, at least over the open subset of $\CP_0$
consisting of line bundles on locally complete intersection subschemes $Z \subset Y$.
As an application we will give an explicit formula for the symplectic form
on the Hilbert scheme of lines on a cubic fourfold.

Let $n = \dim Y$ and $m = n - \deg h$ be the codimension in $Y$ of subschemes
parameterized by $\Hilb^h(Y)$. The tangent space to $\Hilb^h(Y)$
at a point $[Z] \in \Hilb^h(Y)$ is canonically isomorphic
to $H^0(Z,\CN_{Z/Y})$.
Let $\Hilb_{lci}^h(Y) \subset \Hilb^h(Y)$ denote
the open subset of equidimensional locally complete intersection subschemes.
Recall that for any $[Z] \in \Hilb^h_{lci}(Y)$
we have the adjunction exact triple
$$
0 \to \CN^\dual_{Z/Y} \xrightarrow{\kappa} \Omega^1_{Y|Z} \to \Omega_Z \to 0.
$$
Denote by $\kappa \in H^0(Z,\CN_{Z/Y}\otimes\Omega^1_{Y|Z}) = \Hom(\CN^\dual_{Z/Y},\Omega^1_{Y|Z})$
the element corresponding to the first morphism in this triple.
Take any $\omega \in H^{n-m}(Y,\Omega_Y^{n-m+2})$ and denote by
$\omega_{|Z} \in H^{n-m}(Z,\Omega_{Y|Z}^{n-m+2})$ its restriction to $Z$.
We define a $2$-form $\beta_\omega$ on $\Hilb^h_{lci}(Y)$ as follows.
For any $s_1,s_2 \in H^0(Z,\CN_{Z/Y})$ we set
\begin{equation}\label{beta}
\beta_\omega(s_1,s_2) = (\kappa^{\wedge(m-2)}\wedge s_1\wedge s_2\wedge \omega_{|Z}) \ccap [Z],
\end{equation}
where $\kappa^{\wedge(m-2)}$ is the $(m-2)$-fold wedge product of $\kappa$
in $H^0(Z,\wedge^{m-2}\CN_{Z/Y}\otimes\Omega^{m-2}_{Y|Z})$, so that
$\kappa^{\wedge(m-2)}\wedge s_1\wedge s_2\wedge \omega_{|Z} \in
H^{n-m}(Z,\wedge^m\CN_{Z/Y}\otimes\Omega^{n}_{Y|Z}) \cong H^{n-m}(Z,\Omega^m_Z)$.

Let $\CP^\circ \subset p^{-1}(\Hilb^h_{lci}(Y))$ denote the open subset
of $\PPP_0$ consisting of line bundles on subschemes $Z \subset Y$.

\begin{theorem}\label{ab}
We have $\alpha_\omega = p^*\beta_\omega$ on $\CP^\circ$.
\end{theorem}
\begin{proof}
Take any $\CF \in \CP^\circ$ and put $[Z] = p([\CF])$.
By definition of $p$ this means that $\CF \cong i_*F$ for some line bundle $F$
on $Z$, where $i:Z \to Y$ is the embedding and $Z$ is a locally complete intersection in $Y$.
Take any $v_1,v_2 \in \Ext^1(\CF,\CF)$.

By the definition of $\alpha$, we should compute the product
$\At^{m-2}_Y(\CF)\circ v_1\circ v_2$ in $\Ext^m(\CF,\CF\otimes\Omega^{m-2}_Y)$,
then take its trace in $H^m(Y,\Omega^{m-2}_Y)$ and finally couple it with
$\omega \in H^{n-m}(Y,\Omega^{n-m+2}_Y)$.
By Proposition~\ref{exp-fact} the trace factors
through
$$
H^0(Y,\EXT^m(\CF,\CF\otimes\Omega^{m-2}_Y)) \cong
H^0(Z,F^\dual\otimes F\otimes\wedge^m\CN_{Z/Y}\otimes \Omega^{m-2}_{Y|Z}) =
H^0(Z,\wedge^m\CN_{Z/Y}\otimes \Omega^{m-2}_{Y|Z}),
$$
hence it suffices to compute the image of $\At^{m-2}_Y(\CF)\circ v_1\circ v_2$
in $H^0(Z,\wedge^m\CN_{Z/Y}\otimes \Omega^{m-2}_{Y|Z})$
and then apply the canonical map $\mathsf{can}:H^0(Z,\wedge^m\CN_{Z/Y}\otimes \Omega^{m-2}_{Y|Z}) \to H^m(Y,\Omega^{m-2}_Y)$.
By Lemma~\ref{ltgm}
this image coincides with the product of the $(m-2)$-th wedge power of the image of $\At_Y(\CF)$ in
$H^0(Y,\EXT^1(\CF,\CF\otimes\Omega^1_Y)) = H^0(Z,\CN_{Z/Y}\otimes \Omega^1_{Y|Z})$ and of
the images of $v_1$, $v_2$ in $H^0(Y,\EXT^1(\CF,\CF)) = H^0(Z,\CN_{Z/Y})$.
By Theorem~\ref{epsat}~$(iii)$ the image of $\At_Y(\CF)$ in
$H^0(Z,\CN_{Z/Y}\otimes \Omega^1_{Y|Z})$ equals $\kappa$.
The images of $v_i$ are equal to $p_*(v_i)$.
So, by Lemma~\ref{extfg}~$(ii)$ we have
$$
\alpha_\omega(v_1,v_2) =
\mathsf{can}(\kappa^{\wedge(m-2)}\wedge p_*(v_1)\wedge p_*(v_2))\ccup\omega\ccap[Y].
$$
Further, we note that the cap-products $\ccap[Y]$ and $\ccap[Z]$ are nothing but the Serre duality
on $Y$ and $Z$ respectively. Therefore, by Lemma~\ref{can} we have
$$
\mathsf{can}(-)\ccup\omega\ccap[Y] = (-)\ccup\omega_{|Z}\ccap[Z]
$$
which completes the proof.
\end{proof}

Now we apply Theorem~\ref{ab} in the following case.
We take for $Y$ a smooth cubic hypersurface in $\PP^5$ (a cubic 4-fold),
and $h(n) = n + 1$.

\begin{lemma}
If $\CF$ is a semistable sheaf on $Y$ with Hilbert polynomial $h_\CF(n) = n + 1$
then $\CF = \CO_\ell$, the structure sheaf of a line $\ell \subset Y$.
\end{lemma}
\begin{proof}
By Riemann--Roch we have $\mathsf{ch}(\CF) = [\ell]-\frac12[\mathsf{pt}]$,
hence $\CF$ is a rank~1 sheaf on a line. By semistability $\CF$ has
no $0$-dimensional torsion, hence $\CF \cong \CO_\ell(k)$ for some $k \in \ZZ$.
Finally, computing $h_{\CO_\ell(k)}(n) = n + (k+1)$, we conclude that $k = 0$.
\end{proof}

We conclude that $\CP = \Hilb^h(Y) = F(Y)$ is the Hilbert scheme of lines
and the projection map $p:\CP \to \Hilb^h(Y)$ is the identity.
Thus, Theorem~\ref{ab} gives a way to compute the form $\alpha_\omega$.
We are going to do this explicitly.

Recall that $\Omega^4_Y \cong \CO_Y(-3) \cong \CN^\dual_{Y/\PP^5}$.
So, we take the form $\omega \in H^1(Y,\Omega^3_Y)$ corresponding to the
extension $\nu_{Y/\PP^5} \in \Ext^1(\Omega_Y,\CN^\dual_{Y/\PP^5})$
under the isomorphisms
$$
\Ext^1(\Omega_Y,\CN^\dual_{Y/\PP^5}) =
H^1(Y,\CT_Y\otimes\CO_Y(-3)) \cong
H^1(Y,\CT_Y\otimes\Omega^4_Y) \cong
H^1(Y,\Omega^3_Y).
$$
According to~\eqref{beta} we have to compute $\kappa\wedge\omega_{|\ell}$.

\begin{lemma}
The wedge product
$$
\kappa\wedge\omega_{|\ell} \in
H^1(\ell,\CN_{\ell/Y}\otimes\Omega^4_{Y|\ell}) \cong
H^1(\ell,\CN_{\ell/Y}(-3)) \cong
H^1(\ell,\CN_{\ell/Y}\otimes\CN^\dual_{Y/\PP^5|\ell})
$$
is given by the extension class of the normal bundles exact sequence
\begin{equation}\label{norm-b-s}
0 \lra \CN_{\ell/Y} \lra \CN_{\ell/\PP^5} \lra \CN_{Y/\PP^5}|_\ell \lra 0.
\end{equation}
\end{lemma}
\begin{proof}
Consider the following commutative diagram
$$
\xymatrix@R=15pt{
& 0 \ar[d] & 0 \ar[d] \\
& \CN^\dual_{Y/\PP^5}|_\ell \ar@{=}[r] \ar[d] & \CN^\dual_{Y/\PP^5}|_\ell \ar[d] \\
0 \ar[r] & \CN^\dual_{\ell/\PP^5} \ar[r] \ar[d] & \Omega_{\PP^5|\ell} \ar[r] \ar[d] & \Omega_\ell \ar[r] \ar@{=}[d] & 0 \\
0 \ar[r] & \CN^\dual_{\ell/Y} \ar[r]^-\kappa \ar[d] & \Omega_{Y|\ell} \ar[r] \ar[d] & \Omega_\ell \ar[r] & 0 \\
& 0 & 0
}
$$
Note that the form $\omega_{|\ell}$ by definition is given by the restriction
of the extension $\nu_{Y/\PP^5}$ to $\ell$ which coincides with the middle column of the diagram.
Therefore, $\kappa\wedge\omega_{|\ell}$ equals the extension class of the left column.
It remains to note that the extension classes of dual exact sequences coincide.
\end{proof}

So, to compute $\alpha_\omega = \beta_\omega$, it remains
to identify the element $\sigma \in H^1(\ell,\CN_{\ell/Y}(-3))$
corresponding to the extension class of~\eqref{norm-b-s}.

Recall that according to \cite{CG},
there are two possibilities for the normal bundle $\CN_{\ell/Y}$ of a line $\ell$ in $Y$.
We have either $\CN_{\ell/Y} = \CO \oplus \CO \oplus \CO(1)$, or
$\CN_{\ell/Y} = \CO(-1) \oplus \CO(1) \oplus \CO(1)$.
The lines $\ell$ with the normal bundle of the first type fill an open subset of $F(Y)$
and the lines with the normal bundle of the second type fill a subset of codimension 2 in $F(Y)$.

First of all consider the case of $\CN_{\ell/Y} = \CO \oplus \CO \oplus \CO(1)$.
Denote by $e_1,e_2,e_3$ rational sections of twists of $\CN_{\ell/Y}$ such that
$\CN_{\ell/Y}= \CO e_1\oplus \CO e_2\oplus\CO(1)e_3$.
Then the sequence~\eqref{norm-b-s} twisted by $\OOO(-3)$ takes form
$$
0 \to \CO_\ell(-3)e_1 \oplus \CO_\ell(-3)e_2 \oplus \CO_\ell(-2)e_3 \to \CO_\ell(-2)^{\oplus 4} \to \CO_\ell \to 0.
$$
We see that $\sigma$ is a generator of the kernel of the map
$$
H^1(\ell,\CO_\ell(-3)e_1 \oplus \CO_\ell(-3)e_2 \oplus \CO_\ell(-2)e_3) \to H^1(\ell,\CO_\ell(-2)^{\oplus 4}).
$$
This implies that the component $\sigma_3$ of $\sigma$ in the third summand $\CO_\ell(-2)$ of $\CN_{\ell/Y}(-3)$ is zero,
while the components $\sigma_1$ and $\sigma_2$ in the first two summands are linearly independent elements of
$H^1(\ell,\CO_\ell(-3))$.
Using the Serre duality $H^1(\ell,\CO_\ell(-3)) \cong H^0(\ell,\CO(1))^\dual$,
let us endow $\ell$ with homogeneous coordinates
$t_0,t_1$, dual to $\sigma_1,\sigma_2$. The sections $v_i \in H^0(\ell,\CN_\ell)$
can be written as
$$
v_i = a_ie_1 + b_ie_2 + (c_it_0 + d_it_1)e_3.
$$
Then it is clear that
\begin{equation}\label{2-form-in-coo}
\alpha_\omega(v_1,v_2) =
\sigma\wedge v_1\wedge v_2 =
b_1c_2 - b_2c_1 - a_1d_2 + a_2d_1.
\end{equation}

The case of $\CN_{\ell/Y} = \CO(-1) \oplus \CO(1) \oplus \CO(1)$ is considered similarly.
Denote by $e_1,e_2,e_3$ rational sections of twists of $\CN_{\ell/Y}$ such that
$\CN_{\ell/Y}= \CO(-1) e_1\oplus \CO(1) e_2\oplus\CO(1)e_3$.
Then the sequence~\eqref{norm-b-s} twisted by $\CO(-3)$ takes form
$$
0 \to \CO_\ell(-4)e_1 \oplus \CO_\ell(-2)e_2 \oplus \CO_\ell(-2)e_3 \to \CO_\ell(-2)^{\oplus 4} \to \CO_\ell \to 0.
$$
Then $\sigma$ is a generator of the kernel of the map
$$
H^1(\ell,\CO_\ell(-4)e_1 \oplus \CO_\ell(-2)e_2 \oplus \CO_\ell(-2)e_3) \to H^1(\ell,\CO_\ell(-2)^{\oplus 4}).
$$
This implies that the components $\sigma_2$ and $\sigma_3$ of $\sigma$ are zero,
while the component $\sigma_1$ is a nondegenerate element of
$H^1(\ell,\CO_\ell(-4)) \cong S^2H^1(\ell,\CO_\ell(-3))$.
Using the Serre duality $H^1(\ell,\CO_\ell(-3)) \cong H^0(\ell,\CO(1))^\dual$,
let us endow $\ell$ with homogeneous coordinates
$t_0,t_1$, which are isotropic for $\sigma_1$.
The sections $v_i \in H^0(\ell,\CN_\ell)$ can be written as
$$
v_i = (a_it_0 + b_it_1)e_2 + (c_it_0 + d_it_1)e_3.
$$
Then it is clear that
\begin{equation}\label{2-form-in-coo2}
\alpha_\omega(v_1,v_2) =
\sigma\wedge v_1\wedge v_2 =
b_1c_2 - b_2c_1 + a_1d_2 - a_2d_1.
\end{equation}

\section{A 10-dimensional example}\label{sdim10}

Let $X$ be a smooth 3-dimensional cubic hypersurface in
$\PP^4$. By a normal elliptic quintic in $X$, we
mean a curve $C\subset\PP^4$, contained in $X$ and
projectively equivalent to
the image of an elliptic curve $E$ by the linear system
$|5o|$, where $o\in E$ is a point of $E$.
Equivalently, $C$ is a smooth connected curve in $X$ of degree 5
and of genus 1 such that its linear span $\langle C\rangle$
is $\PP^4$, see \cite{H}.

To each $C$ as above, we associate
the vector bundle $\EEE =\EEE_C$ obtained from
$C$ by Serre's construction:
\begin{equation}\label{serre}
0\lra \OOO_X\xrightarrow{\ \ s\ \ } \EEE (1) \xrightarrow{\ \ t\ \ } \III_C(2) \lra 0\; ,
\end{equation}
where $\III_C=\III_{C,X}$ is the ideal sheaf of $C$ in $X$,
and $s$ is a section of $\EEE (1)$ which has $C$ is its
zero locus.
Since the class of $C$ modulo algebraic equivalence
is $5[\ell ]$, where $\ell$ is a line in $X$ and $[\ell ]$ as its class in
$H^4(X,\ZZ )$, (\ref{serre}) implies that
$c_1(\EEE )=0, c_2(\EEE )=2[\ell ]$. Further, $\det\EEE$ is trivial,
and hence  $\EEE$ is self-dual as soon as it is a vector
bundle (that is, $\EEE^\dual\simeq\EEE$).
See \cite[Sect. 2]{MT2} for further details on this construction.
As follows from \cite{IM}, \cite{B1}, \cite{Dr} (see also
\cite{B2}, where the relevant results of the other three
papers are summarized) or \cite{Ku1},
the vector bundles $\EEE$ of this type have several other
equivalent characterizations:

\begin{theorem}\label{bundles-e}
Let $\EEE$ be a rank-$2$ vector bundle on $X$.
Then the following properties are equivalent:

$(i)$ $\EEE$ is stable with Chern classes $c_1=0$, $c_2=2[\ell ]$.

$(ii)$ $\EEE$ is isomorphic to a vector bundle obtained
by Serre's construction \eqref{serre} from
a normal elliptic quintic $C\subset X$.

$(iii)$ $\EEE$ has Chern classes $c_1=0$, $c_2=2[\ell ]$
and the intermediate cohomology of the twists of $\EEE$
vanishes:
$$
H^i(X,\EEE (j))=0\ \mbox{\em for}\ i=1,2\ \mbox{\em and for all}\
j\in\ZZ .
$$

$(iv)$ There exists a Pfaffian representation of $X$, that is
a skew-symmetric $6$ by $6$ matrix $M$ of linear forms on $\PP^4$
such that the equation of $X$ is $\Pf (M)=0$, and $\EEE\simeq \KKK (1)$,
where $\KKK$ is the kernel bundle of $M$: it is defined as a
rank $2$ subbundle of the trivial rank-$6$ bundle $\CO_X^{\oplus 6}$
over $X$ whose fiber $K_x$ over $x\in X$ is the kernel
of the rank-$4$ linear map $M(x):\CC^6\lra\CC^6$. Equivalently,
$\EEE (1)\simeq \CCC$, where $\CCC$ is the cokernel bundle of $M$.

$(v)$ There exists
a skew-symmetric $6$ by $6$ matrix $M$ of linear forms on $\PP^4$
such that $\EEE (1)$ considered as a sheaf on $\PP^4$ can be included
in the following exact sequence:
$$
0\lra \CO(-1)^{\oplus 6} \xrightarrow{M}  \CO^{\oplus 6} \lra\EEE (1)\lra 0.
$$
\end{theorem}

The vector bundles as in the above theorem
possess the following property:

\begin{lemma}\label{hi-ee}
Let $\EEE$ be a vector bundle on $X$ satisfying
any of the equivalent conditions $(i)$--$(v)$
of Theorem~$\ref{bundles-e}$.
Then $H^\bullet(X,\CE) = H^\bullet(X,\CE(-1)) = H^\bullet(X,\CE(-2)) = 0$.
\end{lemma}

\begin{proof}
Follows immediately from Theorem~\ref{bundles-e}~(v).
\end{proof}

Let $M_X=M_X(2;0,2)$ be the moduli space of stable rank-2 vector bundles
with Chern classes $c_1(\EEE )=0, c_2(\EEE )=2[\ell ]$.
There is a natural map
$\phi_X :M\lra J(X)$ to the intermediate Jacobian $J(X)$ of $X$,
well-defined modulo a constant
translation in $J(X)$. It can be described as follows.
According to \cite{Mur}, the Chow group $A_1(0)_X$
of algebraic 1-cycles of degree $0$ on a smooth cubic threefold
$X$ modulo rational equivalence is canonically
isomorphic to the intermediate Jacobian $J(X)$.
Taking any 1-cycle $Z_0$ of degree $d$ on $X$, we obtain also
an identification $A_1(d)_X\lra A_1(0)_X=J(X)$
for the set $A_1(d)_X$ of rational equivalence classes of
degree-$d$ cycles,
$(Z)\in A_1(d)_X\mapsto (Z-Z_0)\in J(X)$,
where $(Z)$ denotes the class of $Z$ modulo rational equivalence.
This is nothing but
the Abel--Jacobi map on the algebraic cycles of degree $d$.

Grothendieck defined in \cite{G} the Chern classes with values
in the Chow groups of algebraic cycles modulo rational equivalence.
Let us denote these Grothendieck--Chern classes by $\fc_i(\EEE )$.
Then the wanted map $\phi_X$ can be
defined by $[\EEE ]\in M_X\mapsto \fc_2(\EEE )-(Z_0)$ for some
fixed reference 1-cycle $Z_0$ of degree 2. One can choose,
for example, $Z_0=2\ell$. Remark that if $\EEE =\EEE_C$ is obtained by
Serre's construction from a normal elliptic quintic $C$,
then $\fc_2(\EEE )= (C) - h^2$, where $h^2$ is the class
of plane cubic curve, a linear section $\PP^2\cap X$.

It follows from the results of \cite{MT1}, \cite{IM} and \cite{Dr} that
$\phi_X$ is an open immersion, thus $M_X(2;0,2)$ is isomorphic
to an open subset of $J(X)$ (see also \cite{B2} and \cite{Ku1}).

Now we are passing to dimension 4.
Let $Y\subset\PP^5$ be a nonsingular cubic fourfold.
Denote by $\FMY$ the moduli space of sheaves on $Y$ of the form $i_*\CE$,
where $[\CE] \in M_X(2;0,2)$, $X$ is a nonsingular hyperplane section of $Y$,
and $i:X \to Y$ is the embedding.
There is a natural map $\pi :\FMY \lra\check{\PP}^5$ whose image
is the complement of $Y^\dual$, the projectively dual variety of~$Y$.
According to \cite{MT2}, $\FMY$ is a nonsingular 10-dimensional variety.

Let $\nu$ be the generator of $H^1(Y, \Omega_Y^3)\simeq\CC$
defined in Theorem~\ref{sf1}, and $\alp_\nu$ the associated
2-form on  $\FMY$. We have already proved its closedness in Section~\ref{c2f}.
Now we will see its nondegeneracy.

\begin{theorem}
The $2$-form $\alpha_\nu$ on the moduli space $\FMY$ is nondegenerate.
\end{theorem}
\begin{proof}
By Lemma~\ref{hi-ee}, we have $\CE \in \CCC_Y$ for every $[\CE] \in \FMY$.
So, Theorem~\ref{sf1} applies.
\end{proof}

\begin{remark} The paper \cite{MT2} provides a
construction of a nondegenerate 2-form
on  $\FMY$, but does not prove its closedness.
It is just the Yoneda pairing $\Lambda$,
as defined by (\ref{yo}), and one can treat it
as a global 2-form on $\FMY$ because
the $1$-dimensional
vector spaces $\Ext^2(i_*\EEE ,i_*\EEE )$ fit into
a trivial line bundle on $\FMY$  as $\EEE$ runs over
$\FMY$. It is also proved in loc. cit.
that $\pi$ is a Lagrangian fibration for $\Lambda$. As
$\alp_\nu$ factors through $\Lambda$, the same holds for $\alp_\nu$.
\end{remark}

\begin{remark}
The second Chern class mappings $\EEE\mapsto \fc_2(\EEE)\in A_1(2)_X$ over the
smooth hyperplane sections $X$ of $Y$
identify $\FMY$ with an open subset of the family $\AAA$ of varieties $A_1(2)_X$. The latter family is an
algebraic torsor
under the relative intermediate Jacobian $\JJJ$ of the family of
smooth hyperplane sections of $Y$. By \cite{DM}, 8.5.2 ,
$\JJJ$ has a natural
symplectic structure $\alp_\JJJ$ such that the map $\JJJ\lra\check{\PP}^5$ is a Lagrangian fibration; let us say for short that
$\alp_\JJJ$ is a Lagrangian structure on $\JJJ/\check{\PP}^5$.
It is easy to see that a Lagrangian structure on a family of abelian varieties
induces a Lagrangian structure on any its algebraic torsor.
Let us denote the thus induced Lagrangian structure on $\AAA/\check{\PP}^5$
by $\alp_\AAA$. Then it is plausible that $\alp_\nu$ coincides
with the restriction of $\alp_\AAA$ up to a constant factor.
A way to prove this might be to find a partial compactification $\bar\FMY$
of $\FMY$, such that $h^0(\Omega^2_{\bar\FMY})=1$ and both $\alp_\nu$, $\alp_\AAA$ extend to $\bar\FMY$.
\end{remark}


\begin{thebibliography}{MMMMM}


\bibitem[AB]{AB} Atiyah, M. F., Bott, R.: {\em The Yang-Mills equations over Riemann surfaces},  Philos. Trans. Roy. Soc. London Ser. {\bf A  308}  (1983), 523--615.

\bibitem[AK]{AK} Altman, A. B., Kleiman, S. L.:
{\em Compactifying the Picard scheme,}
Adv. in Math. {\bf 35} (1980),  50--112.


\bibitem[ALJ]{ALJ} Ang\'eniol, B., Lejeune-Jalabert, M.: {\em
Calcul diff\'erentiel et classes caract\'eristiques en
g\'eom\'etrie alg\'ebrique},
Travaux en Cours, 38, Hermann, Paris, 1989.

\bibitem[Ar]{Ar}   Artamkin, I. V.: {\em On the deformation of sheaves},
Math. USSR Izv. {\bf 32} (1989), 663–668.

\bibitem[At]{At} Atiyah, M. F.: {\em
Complex analytic connections in fiber bundles,}
Trans. Amer. Math. Soc. {\bf 85} (1957), 181--207.


\bibitem[B1]{B1} Beauville, A.: {\em
Determinantal hypersurfaces,} Michigan Math. J.
{\bf 48} (2000), 39--64.

\bibitem[B2]{B2} Beauville, A.: {\em
Vector bundles on the cubic threefold}, Bertram, Aaron (ed.) et al.,
Symposium in honor of C. H. Clemens,
University of Utah, March 2000, Contemp. Math. {\bf 312}, 71--86 (2002).



\bibitem[BD]{BD} Beauville, A., Donagi, R.: {\em La vari\'et\'e des
droites d'une hypersurface cubique de dimension 4,} C. R. Acad. Sci. Paris
Ser. I Math. {\bf 301}, no. 14, 703--706 (1985).

\bibitem[Bei]{Bei} Beilinson, A.:
{\em Coherent sheaves on $P\sp{n}$ and problems in linear algebra} (Russian)
Funktsional. Anal. i Prilozhen. {\bf 12} (1978), no. 3, 68--69.

\bibitem[Bon]{Bon} Bondal, A.:
{\em Representations of associative algebras and coherent sheaves},
(Russian)  Izv. Akad. Nauk SSSR Ser. Mat. {\bf 53} (1989), no. 1, 25--44;
translation in  Math. USSR-Izv. {\bf 34}  (1990),  no. 1, 23--42.

\bibitem[BK]{BK} Bondal, A., Kapranov, M.:
{\em Representable functors, Serre functors, and reconstructions},
(Russian)  Izv. Akad. Nauk SSSR Ser. Mat. {\bf 53} (1989), no. 6, 1183--1205, 1337;
translation in  Math. USSR-Izv. {\bf 35} (1990), no. 3, 519--541.

\bibitem[Bot-1]{Bot-1} Bottacin, F.: {\em
Poisson structures on moduli spaces of sheaves over Poisson surfaces},
Invent. Math. {\bf 121} (1995), 421--436.

\bibitem[Bot-2]{Bot-2} Bottacin, F.: {\em Symplectic geometry on moduli spaces of stable pairs},  Ann. Sci. E'cole Norm. Sup. (4) {\bf 28}  (1995),  391--433.


\bibitem[Bot-3]{Bot-3} Bottacin, F.: {\em
Poisson structures on moduli spaces of parabolic bundles on surfaces}, Manuscripta Math. {\bf  103}  (2000),   31--46.

\bibitem[BR]{BR} Biswas, I., Ramanan, S.: {\em
An infinitesimal study of the moduli of Hitchin pairs},
J. London Math. Soc. (2) {\bf 49} (1994), 219–231.



\bibitem[BuF1]{BuF1}  Buchweitz, R.-O.,  Flenner, H.: {\em
The Atiyah-Chern character yields the semiregularity map
as well as the infinitesimal Abel-Jacobi map}, The
arithmetic and geometry of algebraic cycles (Banff, AB, 1998),
33--46, CRM Proc. Lecture Notes, 24, Amer. Math. Soc., Providence, 2000.

\bibitem[BuF2]{BuF2}  Buchweitz, R.-O.,  Flenner, H.: {\em
A Semiregularity Map for Modules and Applications to Deformations},
Compositio Math. {\bf 137}  (2003), 135--210.


\bibitem[CG]{CG} Clemens, C. H., Griffiths, P. A.: {\em The
intermediate Jacobian of the cubic threefold,} Annals of Math.
{\bf 95}(1972), 281--356.

\bibitem[dJS]{dJS} de Jong, A. J.; Starr, J.: {\em Cubic fourfolds and spaces of rational curves},  Illinois J. Math. {\bf 48}  (2004),  415--450.

\bibitem[DM]{DM} Donagi, R., Markman, E.: {\em Spectral covers, algebraically
completely integrable Hamiltonian systems, and moduli of bundles,}
In: Francaviglia, M. (ed.) et al., Integrable systems and quantum groups,
CIME Lectures, Italy, June 14-22, 1993. Lect. Notes Math., 1620, Springer-Verlag,
Berlin, 1996, 1--119.

\bibitem[Dr]{Dr} Druel, S.: {\em
Espace des modules des faisceaux semi-stables de rang 2 et
de classes de Chern $c_{1}=0$, $c_{2}=2$ et $c_{3}=0$ sur une
hypersurface cubique lisse de $\mathbb{P}^{4}$},
Internat. Math. Res. Notices {\bf 2000}, No. 19, 985--1004.

\bibitem[Gri]{Gri} Griffiths, P. A.: {\em
Periods of integrals on algebraic manifolds, II, Local study of the period mapping},  Amer. J. Math. {\bf 90}  (1968), 805--865.


\bibitem[GR]{GR} Gorodentsev, A., Rudakov, A.:
{\em Exceptional vector bundles on projective spaces},
Duke Math. J. {\bf 54} (1987), no.~1, 115--130.



\bibitem[Gro]{G} Grothendieck, A.: {\em
La th\'eorie des classes de Chern},
Bull. Soc. Math. France {\bf 86} (1958), 137--154.

\bibitem[Ha]{Ha-1} Hartshorne, R.: {\em
Residues and duality,}
A seminar on the work of A. Grothendieck, Harvard, 1963/64,
Lecture Notes in Math., {\bf 20} Springer, Berlin--N. Y., 1966.


\bibitem[Hi]{Hi} Hitchin, N.: {\em Stable bundles and integrable systems,}
Duke Math. J.  {\bf 54}  (1987),  91--114.

\bibitem[Hu]{H} Hulek, K.: {\em Projective geometry of elliptic
curves,} Ast\'erisque {\bf 137} (1986).


\bibitem[HL]{HL} Huybrechts, D., Lehn, M.: {\em
The geometry of moduli spaces of sheaves,} Aspects of Mathematics,
E31, Braunschweig, Vieweg, 1997.

\bibitem[IIS]{IIS} Inaba, M., Iwasaki, K., Saito, M.-H.: {\em
Moduli of Stable Parabolic Connections, Riemann-Hilbert correspondence and Geometry of Painleve' equation of type VI}, I: math.AG/0309342, II:
math.AG/0605025.

\bibitem[IMan]{IMan} Iliev, A., Manivel, L.: {\em
Cubic hypersurfaces and integrable systems},
math.AG/0605260.

\bibitem[IMar]{IM} Iliev, A., Markushevich, D.: {\em The Abel--Jacobi map
for  a cubic threefold and periods of Fano
threefolds of degree 14}, Doc. Math. {\bf 5} (2000), 23--47.

\bibitem[Ill]{Ill} Illusie, L.: {\em
Complexe cotangent et d\'eformations}, I, II, Lecture Notes in
Math., 239, 283, Springer Verlag, Berlin-Heidelberg- New York, 1971, 1972.

\bibitem[K]{K} Kobayashi, S.: {\em
Simple vector bundles over symplectic Kaehler manifolds},
Proc. Japan Acad., Ser. A {\bf 62} (1986), 21--24.

\bibitem[KN]{KN} Kronheimer, P. B., and Nakajima, H.:
{\em Yang-Mills instantons on ALE gravitational instantons},
Math. Ann. {\bf 288:2} (1990), 263--307.

\bibitem[Ku1]{Ku1} Kuznetsov, A. G.: {\em Derived category of a cubic threefold and the
variety $V\sb {14}$}, Proc. Steklov Inst. Math. {\bf 246} (2004), 171--194.

\bibitem[Ku2]{Ku2} Kuznetsov, A. G.: {\em Homological projective duality for Grassmannians of lines},
math.AG/0610957.

\bibitem[LP]{LP} Le Potier, J.: {\em
Faisceaux semi-stables et syst\`emes coh\'erents},
Vector bundles in algebraic geometry
(Durham, 1993), Ed. N. Hitchin et al., 179--239,
London Math. Soc. Lecture Note Ser., 208,
Cambridge Univ. Press, Cambridge, 1995.








\bibitem[MT1]{MT1} Markushevich, D., Tikhomirov, A. S.:
{\em The Abel--Jacobi map of a moduli component of vector bundles
on the cubic threefold},  J. Algebraic Geom.
{\bf 10} (2001), 37--62.

\bibitem[MT2]{MT2} Markushevich, D., Tikhomirov, A. S.:
{\em Symplectic structure on a moduli space of sheaves on
the cubic fourfold},
Izv. Ross. Akad. Nauk Ser. Mat. {\bf 67} (2003), 131--158.

\bibitem[Muk-1]{Muk-1} Mukai, S.: {\it Symplectic structure of the moduli
of sheaves on an abelian or $K3$ surface,} Invent. Math. {\bf 77} (1984),
101--116.

\bibitem[Muk-2]{Muk-2} Mukai, S.: {\it Moduli of vector bundles on $K3$ surfaces and symplectic manifolds},
Sugaku Expositions {\bf 1} (1988), 139--174.

\bibitem[Mur]{Mur} Murre, J. P.: {\em
Some results on cubic threefolds}, In: Classification of algebraic
varieties and compact complex manifolds, 140--160,
Lecture Notes in Math., Vol. 412, Springer, Berlin, 1974.


\bibitem[O'G-1]{O'G} O'Grady, K. G.: {\em
Algebro-geometric analogues of Donaldson's polynomials},
Invent. Math. {\bf 107} (1992), 351--395.

\bibitem[O'G-2]{OG} O'Grady, K. G.: {\em Desingularized moduli spaces of sheaves on a $K3$},  J. Reine Angew. Math. {\bf 512}  (1999), 49--117.


\bibitem[OSS]{OSS} Okonek, C., Schneider, M., Spindler, H.:
{\em Vector bundles on complex projective spaces}
Progress in Mathematics, 3.
Birkha\"user, Boston, Mass., 1980. vii+389 pp.

\bibitem[Ran]{Ran} Ran, Z.: {\em On the local geometry of moduli spaces of locally free sheaves}, In:  Moduli of vector bundles (Sanda, 1994; Kyoto, 1994),  213--219, Lecture Notes in Pure and Appl. Math., 179, Dekker, New York, 1996.



\bibitem[Si]{S} Simpson, C. T.: {\em Moduli of representations of the
fundamental group of a smooth projective variety I,} Publ. Math.
I.H.E.S. {\bf 79} (1994), 47--129.

\bibitem[T]{T} Thomas, R. P.: {\em
A holomorphic Casson invariant for Calabi-Yau 3-folds, and bundles on $K3$ fibrations},
J. Differential Geom. {\bf 54} (2000),  367--438.



\bibitem[Tyu]{Tyu} Tyurin, A.N.: {\em
Symplectic structures on the varieties of moduli of vector
bundles on algebraic surfaces with
$p_g>0$ },
Math. USSR, Izv. {\bf 33} (1989), 139--177 $=$ Izv. Akad. Nauk SSSR, Ser. Mat.
{\bf 52} (1988), 813--852.

\bibitem[Ve]{Ve} Verdier, J.-L.:
{\em Des cat\'egories d\'eriv\'ees des cat\'egories ab\'eliennes},
Ast\'erisque  No. 239  (1996), xii+253 pp.


\bibitem[W]{W} Welters, G.E.: {\em Abel--Jacobi isogenies for
certain types of Fano threefolds,} Mathematical Centre Tracts {\bf 141},
Amsterdam, 1981.


\end{thebibliography}
\end{document}